\documentclass[11pt]{amsart}

\calclayout

\usepackage{color}
\usepackage{verbatim}
\usepackage[noadjust]{cite}

\usepackage[utf8x]{inputenc}
\usepackage{t1enc}
\usepackage[english]{babel}
\usepackage{pgfplots}
\pgfplotsset{compat=1.15}
\usepackage{mathrsfs}
\usetikzlibrary{arrows}

\definecolor{ffqqqq}{rgb}{1.,0.,0.}
\definecolor{uuuuuu}{rgb}{0.26666666666666666,0.26666666666666666,0.26666666666666666}

\usepackage{amsmath,amssymb,amsthm}
\usepackage{mathrsfs,esint}
\usepackage{graphicx,wrapfig}
\usepackage[format=hang]{caption}
\newcommand{\R}{\mathbb{R}}

\newcommand{\N}{\mathbb{N}}
\newcommand{\Z}{\mathbb{Z}}

\newcommand{\B}{\mathcal{B}}

\newcommand{\Ha}{\mathcal{H}}

\newcommand{\Leb}{\mathcal{L}}

\newcommand{\eps}{\varepsilon}

\newcommand{\loc}{\text{loc}}

\newcommand{\phii}{\varphi}

\newcommand{\bmat}{\begin{bmatrix}}
\newcommand{\emat}{\end{bmatrix}}

\newcommand{\wtil}{\widetilde}

\newcommand{\st}{\text{ s.t. }}

\providecommand*{\vint}[1]{\mathchoice
          {\mathop{\vrule width 5pt height 3 pt depth -2.5pt
                  \kern -9pt \kern 1pt\intop}\nolimits_{\kern -5pt{#1}}}
          {\mathop{\vrule width 5pt height 3 pt depth -2.6pt
                  \kern -6pt \intop}\nolimits_{\kern -3pt{#1}}}
          {\mathop{\vrule width 5pt height 3 pt depth -2.6pt
                  \kern -6pt \intop}\nolimits_{\kern -3pt{#1}}}
          {\mathop{\vrule width 5pt height 3 pt depth -2.6pt
                  \kern -6pt \intop}\nolimits_{\kern -3pt{#1}}}}

\DeclareMathOperator{\dist}{dist}
\DeclareMathOperator{\diam}{diam}
\DeclareMathOperator{\rad}{rad}

\DeclareMathOperator{\Id}{Id}

\numberwithin{equation}{section}
\theoremstyle{plain}
\newtheorem{thm}[equation]{Theorem}

\newtheorem{prop}[equation]{Proposition}
\newtheorem{cor}[equation]{Corollary}
\newtheorem{lem}[equation]{Lemma}

\theoremstyle{definition}

\newtheorem{defn}[equation]{Definition}

\newtheorem{remark}[equation]{Remark}

\makeatletter
\def\blfootnote{\xdef\@thefnmark{}\@footnotetext}
\makeatother

\begin{document}

\title[Besov trace and extension theorems]{Trace and extension theorems for Besov functions in doubling metric measure spaces}

\blfootnote{2020 {\it Mathematics Subject Classification.} Primary: 46E36 Secondary: 46E35, 30L15}
\blfootnote{{\it Keywords and phrases.} Besov space, trace, extension, hyperbolic filling, metric measure space.}
\blfootnote{{\it Acknowledgments.} The research of I.C.\ was partially supported by the grant PID2022-138758NB-I00 (Spain), as well as the Institute of Mathematics of the Polish Academy of Sciences, and the research of J.K.\ was partially supported by the NSF Grant DMS-\#2054960.  This project was initiated when the first author visited the second at the University of Cincinnati during the spring of 2023, and he would like to thank UC for its warm welcome to a great research environment, where many unforgettable moments and mathematical discussions were shared.  The authors would like to thank Nageswari Shanmugalingam for valuable discussions and suggestions regarding this project, as well as Lizaveta Ihnatsyeva for pointing out Lemma~\ref{lem:Sums}.}

\author{Iv\'an Caama\~no, Josh Kline}

\date{October 1, 2025}

\maketitle

 \begin{abstract}  
 In the setting of a non-complete doubling metric measure space $(\Omega,d,\mu)$, we construct various bounded linear trace and extension operators for homogeneous and inhomogeneous Besov spaces $B^\alpha_{p,q}$. Equipping the boundary $\partial\Omega:=\overline\Omega\setminus\Omega$ with a measure which is codimension $\theta$ Ahlfors regular with respect to $\mu$, these operators take the form
 \[
 T:B^\alpha_{p,q}(\Omega)\to B^{\alpha-\theta/p}_{p,q}(\partial\Omega),\quad E:B^\alpha_{p,q}(\partial\Omega)\to B^{\alpha+\theta/p}_{p,q}(\Omega).
 \]
The trace operators are first constructed under the additional assumption that $\Omega$ is a uniform domain in its completion.  We then use such results along with the technique of hyperbolic filling to remove this assumption in the case that $\Omega$ is bounded. This extends to the doubling setting some earlier results of Marcos and Saksman-Soto proven under the assumption that the ambient measure is Ahlfors regular.    
\end{abstract}

\section{Introduction}
In the setting of a  doubling metric measure space $(X,d,\mu)$, we study various trace and extension results for the Besov spaces $B^\alpha_{p,q}(X,\mu)$, defined by finiteness of the following seminorm:
\[
\Vert u\Vert_{HB^\alpha_{p,q}(X,\mu )}:= \left(\int_0^\infty\left(\int_X\fint_{B(x,t)}|u(y)-u(x)|^pd\mu(y)d\mu(x)\right)^{q/p}\frac{dt}{t^{\alpha q+1}}\right)^{1/q},\]
for $u\in L^1_\loc(X)$, $0<\alpha<1$, $1\le p<\infty$, and $1\le q\le\infty$.  See Section~\ref{section:Besov} for the precise definitions. These nonlocal spaces are closely connected with the study of the fractional Laplacian, nonlocal minimal surfaces, fractional Hardy inequalities, as well as the study of jump processes in connection with Dirichlet forms, for a sampling see \cite{CS,CSt,CRS,D, IMV,BBCK,C,HiKu}.    

Relevant to the topic of this paper, the Besov spaces   arise naturally as the trace class of Sobolev spaces for sufficiently regular domains in $\R^n$, when $p=q$.  For example, it was shown by Gagliardo in \cite{G} that there exist bounded linear trace and extension operators between $W^{1,p}(\Omega)$ and $B^{1-1/p}_{p,p}(\partial\Omega,\Ha^{n-1})$ when $\Omega\subset\R^n$ is a Lipschitz domain.  Furthermore, when a closed set $F\subset\R^n$ is a \emph{$d$-set}, that is, $\Ha^d$ is Ahlfors $d$-regular on $F$, then restrictions of functions in $W^{1,p}(\R^n)$ to $F$ belong to $B^{1-(n-d)/p}_{p,p}(F,\Ha^d)$, as shown by Jonsson-Wallin \cite{JW}.  Thus, the dimension of the boundary set determines the regularity of the trace class.  

In recent decades, a rich theory of analysis on metric measure spaces has been developed, in particular regarding Sobolev space theory in the nonsmooth setting.  One natural substitute for the space $W^{1,p}$ is the Newton-Sobolev class $N^{1,p}(X)$ introduced in \cite{S}.  When $(X,d,\mu)$ is doubling and supports a Poincar\'e inequality, the Newton-Sobolev class provides a fruitful framework from which to study potential theory, functions of bounded variation, quasiconformal maps, and analogs of PDEs, for example, in the metric setting.  For more on these spaces, see the monographs \cite{H,HKST,BB} for example.    

In relation to the study of energy minimization and boundary value problems in this context, the trace theory for the Newton-Sobolev class has undergone significant study in recent years.  It was shown by Mal\'y in \cite{M} that when $(X,d,\mu)$ is doubling and supports a $p$-Poincar\'e inequality, and when $\Omega\subset X$ is a uniform domain whose boundary is bounded, then there exist bounded linear trace and extension operators between $N^{1,p}(\Omega,\mu)$ and $B^{1-\theta/p}_{p,p}(\partial\Omega,\nu)$.  Here, $\partial\Omega$ is equipped with a measure $\nu$ which is codimension $\theta$ Ahlfors regular with respect to $\mu$ for some $\theta>0$.  When $(X,d,\mu)$ is a doubling metric measure space, this codimensionality condition, see Definition~\ref{def:CDAR}, serves as a replacement for the notion of a $d$-set, as used in the Euclidean setting in \cite{JW}.  Furthermore, it was shown by Bj\"orn-Bj\"orn-Shanmugalingam in \cite{BBS} that given any compact doubling metric measure space $(Z,d_Z,\nu)$ and $0<\alpha<1$, there exists a doubling metric measure space $(X,d,\mu)$ supporting a $1$-Poincar\'e inequality such that $X$ is a uniform domain in its completion, $Z$ is bi-Lipschitz equivalent to $\partial X$, and $B^\alpha_{p,p}(Z,\nu)$ is the trace class for $N^{1,p}(X,\mu)$.  The space $(X,d,\mu)$ is a so-called \emph{uniformized hyperbolic filling} of $Z$, see Section~\ref{sec:hyp fill}, and is constructed in such a manner so that $\nu$ is codimension $p(1-\alpha)$ Ahlfors regular with respect to $\mu$.  Recently, \cite{Nages-Ryan} extended these trace and extension results of \cite{M} between the Newton-Sobolev and Besov spaces to the case when the $\Omega$ is an unbounded uniform domain, equipped with a doubling measure supporting a Poincar\'e inequality.

In this paper, we are interested in trace and extension results for Besov spaces themselves, rather than between Newton-Sobolev spaces and Besov spaces. As we consider only these nonlocal function spaces, we do not require that our space supports a Poincar\'e inequality, only assuming that the measure is doubling.  In the Euclidean setting, such results between Besov spaces were also established by Jonsson and Wallin in \cite{JW}, where they proved restriction and extension theorems between $B^\alpha_{p,q}(\R^n)$ and $B^{\alpha-(n-d)/p}_{p,q}(F,\Ha^d)$ for $d$-sets $F\subset\R^n$.  In Ahlfors $Q$-regular metric measure spaces $(X,d,\mu)$, where for all balls $B(x,r)\subset X$, one has 
\[
\mu(B(x,r))\simeq r^Q,
\]
analogous results were established by Marcos \cite{Marcos} via interpolation techniques, and by Saksman-Soto \cite{SS}, using an equivalent formulation of the Besov seminorm defined through hyperbolic fillings.  The trace results in these papers take the form of restriction theorems between globally defined Besov functions and those defined on $d$-sets.

A primary goal of this paper is to construct trace operators on Besov spaces in metric measure spaces whose measure is doubling rather than Ahlfors regular, weakening the assumptions of \cite{Marcos,SS}. A particular motivation for this investigation is the recent interest in the construction of fractional $p$-Laplace-type operators in the doubling metric measure space setting, see for example \cite{CKKSS,EGKSS,CGKS,CGKS2}, and formulation of Dirichlet-type problems for such operators \cite{KLS}. Since the Besov spaces $B^\alpha_{p,p}(\Omega)$ form the natural domain for such operators, the construction of trace operators on Besov spaces is essential to the study of well-posedness of the corresponding Dirichlet problems. The study of well-posedness in \cite{KLS} required such a trace result for $B^\alpha_{p,p}(\Omega)$, see \cite[Theorem~1.1]{KLS}, which was obtained by leveraging known trace and extension results between Newton-Sobolev and Besov spaces via hyperbolic fillings \cite{BBS}.  In this paper, we also consider the case $p\ne q$, and so we cannot take advantage of this relationship with Newton-Sobolev spaces.    

Mirroring the trace theorems for Newton-Sobolev spaces on domains obtained in \cite{M,Nages-Ryan}, we first establish trace theorems under the additional geometric assumption that the domain is uniform in order to obtain the appropriate energy bounds.  In contrast to the uniform domains considered in \cite{M,Nages-Ryan}, however, we allow for possibly non-rectifiable uniform curves in the definition, and so our approach is slightly different, see Definition~\ref{def:Uniform domain} and Remark~\ref{rem:Uniform Domain}.
 The first of our main results is the following, pertaining to the homogeneous Besov spaces:

\begin{thm}\label{thm:Trace}
   Let $(\Omega,d,\mu)$ be a locally compact, non-complete metric measure space, with $\mu$ a doubling measure, such that $\Omega$ is a uniform domain in its completion $\overline\Omega$.  Suppose also that $\partial\Omega:=\overline\Omega\setminus\Omega$, the boundary of $\Omega$, is equipped with a Borel measure $\nu$ which is codimension $\theta$ Ahlfors regular with respect to $\mu$ for some $\theta>0$.  Let $0<\alpha<1$,  $1\le p<\infty$, and $1\le q\le\infty$ be such that $0<\alpha-\theta/p<1$.  Then there exists a bounded linear trace operator
    \[
    T:HB^\alpha_{p,q}(\Omega,\mu)\to HB^{\alpha-\theta/p}_{p,q}(\partial\Omega,\nu)
    \]
    such that for all $u\in HB^\alpha_{p,q}(\Omega,\mu)$,   we have
    \begin{equation}\label{eq:intro trace condition}
    \lim_{r\to 0^+}\fint_{B(z,r)}|u-Tu(z)|^pd\mu=0
    \end{equation}
    for $\nu$-a.e.\ $z\in\partial\Omega$, and 
    \begin{equation}\label{eq:hom trace energy bound}
        \|Tu\|_{HB^{\alpha-\theta/p}_{p,q}(\partial\Omega,\nu)}\le C\|u\|_{HB^\alpha_{p,q}(\Omega,\mu)}, 
    \end{equation}
    where the constant $C\ge 1$ depends only on $\alpha$, $p$, $q$, and $\theta$, as well as the doubling, codimensional Ahlfors regularity, and uniform domain constants. 
    \end{thm}

    When $\partial\Omega$ is additionally assumed to be bounded, the operator $T$ obtained in Theorem~\ref{thm:Trace} also acts as a bounded linear trace operator on the inhomogeneous Besov spaces:
    
    \begin{thm}\label{thm:Intro Inhom Trace}
    Under the same assumptions as Theorem~\ref{thm:Trace}, suppose in addition that $\partial\Omega$ is bounded. Then there exists a bounded linear trace operator  
    \[
    T:B^\alpha_{p,q}(\Omega,\mu)\to B^{\alpha-\theta/p}_{p,q}(\partial\Omega,\nu)
    \]
    such that for all $u\in B^\alpha_{p,q}(\Omega,\mu)$,  \eqref{eq:intro trace condition} and \eqref{eq:hom trace energy bound} hold,  and 
    \begin{equation}\label{eq:trace inhom Lp}
        \|Tu\|_{L^p(\partial\Omega,\nu)}\le C\left(\|u\|_{L^p(\Omega,\mu)}+\|u\|_{HB^\alpha_{p,q}(\Omega,\mu)}\right),
    \end{equation}
    where the constant $C\ge 1$ depends only on 
    $\alpha,\,p,\,q,\,\theta$, $\diam(\partial\Omega)$, as well as the doubling and codimensional Ahlfors regularity constants. 
\end{thm}
Theorem~\ref{thm:Trace} and Theorem~\ref{thm:Intro Inhom Trace} are proved in Sections~\ref{section:Trace} via Proposition~\ref{prop:tracelimit} and Proposition~\ref{thm:TraceEnergy}. 

We will later use Theorem~\ref{thm:Trace} and Theorem~\ref{thm:Intro Inhom Trace}, along with corresponding extension results, to establish a trace theorem on bounded spaces \emph{without} the uniform domain assumption by using the hyperbolic filling construction,  see Theorem~\ref{thm:Intro gen trace}.  However, we believe Theorem~\ref{thm:Trace} and Theorem~\ref{thm:Intro Inhom Trace} may be of independent interest, as the operators there are defined intrinsically, without using this auxiliary construction, and these results apply to cases where $\Omega$ may be unbounded.

A second goal of this paper is to construct corresponding extension operators for Besov spaces in doubling metric measure spaces.  In contrast to Theorem~\ref{thm:Trace} and Theorem~\ref{thm:Intro Inhom Trace}, no uniform domain assumption is needed, and as mentioned above these extension results will be used to construct the trace operator in Theorem~\ref{thm:Intro gen trace} without the uniform domain assumption.  Our operators are constructed using standard Whitney extension techniques, involving Whitney covers, Lispchitz partitions of unity, and discrete convolution, all of which are available in the doubling metric measure space setting.  This general strategy was used in \cite{JW} to construct extension operators from $B^{\alpha-(n-d)/p}_{p,q}(F,\Ha^d)$  to $B^\alpha_{p,q}(\R^n)$ for $d$-sets $F\subset\R^n$. There, however, $C^\infty$ functions were used in the partition of unity, and so to obtain the necessary energy estimates, the authors were able to take advantage of smoothness properties of the discrete convolution along line segments joining certain points in the domain.  Such techniques are not available in the metric setting, and so we utilize the Lipschitz partition of unity instead.

 In the doubling metric measure space setting, similar Whitney extension operators, via Lipschitz partitions of unity, were constructed in \cite{Marcos} from $B^{\alpha-\theta/p}_{p,q}(F)$ to $B^\alpha_{p,q}(X)$ for $\theta$-codimensional closed sets $F\subset X$.  A notable difference between our extension results and that of \cite{Marcos} involves the bounds of the extension operator.  In both \cite{Marcos} and the analogous Euclidean results in \cite{JW}, only the inhomogeneous Besov spaces were considered.  As such, the bounds on the Besov energy of the extension operator were given in terms of both the Besov energy \emph{and} the $L^p$-norm of the original function.  The proof of this bound is included in the preprint version of \cite{Marcos}, though omitted from the published version.  We first construct an extension operator for the homogeneous Besov spaces, in which we bound the energy of the extension by the energy of the original function alone:

\begin{thm}\label{thm:Extension}
Let $(\Omega,d,\mu)$ be a locally compact, non-complete metric measure space, with $\mu$ a doubling measure.  Suppose also that $\partial\Omega:=\overline\Omega\setminus\Omega$, the boundary of $\Omega$, is equipped with a Borel measure $\nu$ which is codimension $\theta$ Ahlfors regular with respect to $\mu$ for some $\theta>0$.  Let $0<\alpha<1$, $1\leq p <\infty$, and $1\le q\le\infty$ be such that $\alpha<1-\theta/p$.  Then there exists a bounded linear extension operator 
\[
\wtil E:HB^\alpha_{p,q}(\partial\Omega,\nu)\to HB^{\alpha+\theta/p}_{p,q}(\Omega,\mu).
\] 
In particular, there exists $C\ge 1$, depending only on $\alpha$, $p$, $q$, $\theta$, and the doubling and codimensional Ahlfors regularity constants such that 
\begin{equation}\label{eq:hom extension energy bound}
\|\wtil Ef\|_{HB^{\alpha+\theta /p}_{p,q}(\Omega,\mu)}\le C\|f\|_{HB^\alpha_{p,q}(\partial\Omega,\nu)},
\end{equation}
for all $f\in HB^\alpha_{p,q}(\partial\Omega,\nu)$.  
\end{thm}

When $\Omega$ is bounded, the extension operator $\wtil E$ constructed for the homogeneous Besov spaces becomes an extension operator for the inhomogeneous spaces:

\begin{cor}\label{cor:inhom ext bounded domain}
    Under the assumptions of Theorem~\ref{thm:Extension}, suppose in addition that $\Omega$ is bounded. Then there exists a bounded linear extension operator 
\[
\wtil E:B^\alpha_{p,q}(\partial\Omega,\nu)\to B^{\alpha+\theta/p}_{p,q}(\Omega,\mu)
\]
such that for all $f\in B^\alpha_{p,q}(\partial\Omega,\nu)$, \eqref{eq:hom extension energy bound} holds and 
\begin{equation}\label{eq:inhom Lp bounds}
    \|\wtil Ef\|_{L^p(\Omega,\mu)}\le C\|f\|_{L^p(\partial\Omega,\nu)},
\end{equation}
where $C\ge 1$ depends only on $p$, $\theta$, $\diam(\Omega)$, as well as the doubling and codimensional Ahlfors regularity constants.
\end{cor}

If $\Omega$ is not assumed to be bounded, then by applying a Lipschitz cutoff function to the homogeneous extension above, we construct the following extension operator for the inhomogeneous spaces, analogous to that obtained in \cite{Marcos}, where the Besov energy of the extension is controlled by both the Besov energy and $L^p$-norm of the original function.  Our proof is similar to the Euclidean proof in \cite{JW}, where a smooth cutoff function was used.

\begin{thm}\label{thm:InhomExtension}
  Under the same assumptions as Theorem~\ref{thm:Extension}, there exists a bounded linear extension operator 
\[
 E:B^\alpha_{p,q}(\partial\Omega,\nu)\to B^{\alpha+\theta/p}_{p,q}(\Omega,\mu).
\] 
In particular, there exists $C\ge 1$, depending only on $\alpha$, $p$, $q$, $\theta$,  as well as the doubling and codimensional Ahlfors regularity constants, such that 
\begin{align*}
\|Ef\|_{L^p(\Omega,\mu)}&\le C\|f\|_{L^p(\partial\Omega,\nu)}\quad\text{ and }\\
\|Ef\|_{HB^{\alpha+\theta/p}_{p,q}(\Omega,\mu)}&\le C\left(\|f\|_{HB^\alpha_{p,q}(\partial\Omega,\nu)}+\|f\|_{L^p(\partial\Omega,\nu)}\right)
\end{align*}
for all $f\in B^\alpha_{p,q}(\partial\Omega,\nu).$  
\end{thm}

In Corollary~\ref{cor:Identity}, we show that under the assumptions of Theorem~\ref{thm:Trace}, we have that $T\circ E$ and $T\circ\wtil E$ act as identity operators on the Besov spaces on the boundary.

Having established the relevant extension operators, we then obtain the following trace result without the uniform domain assumption:

\begin{thm}\label{thm:Intro gen trace}
   Let $(Z,d,\nu)$ be a locally compact, non-complete, bounded metric measure space,with $\nu$ a doubling measure.  Suppose that $\partial Z:=\overline Z\setminus Z$, the boundary of $Z$, is equipped with a Borel measure $\pi$, which is codimension $\theta$ Ahlfors regular with respect to $\nu$ for some $\theta >0$. Let $1\le p<\infty$, $1\le q\le\infty$ and $\theta/p<\alpha<1$.   Then there exist a bounded, linear trace operator
    \[
    T:B^\alpha_{p,q}(Z,\nu)\to B^{\alpha-\theta/p}_{p,q}(\partial Z,\pi),
    \]
    such that for all $u\in B^\alpha_{p,q}(Z,\nu)$, we have
    \begin{equation}\label{eq:Intro gen Trace condition}
        \lim_{r\to 0^+}\fint_{B(z,r)}|u-Tu(z)|^pd\nu=0
    \end{equation}
    for $\pi$-a.e.\ $z\in \partial Z$, and
    
    \begin{equation*}
        \|Tu\|_{HB^{\alpha-\theta/p}_{p,q}(\partial Z,\pi)}\le C\|u\|_{HB^\alpha_{p,q}(Z,\nu)},\qquad\|Tu\|_{L^p(\partial Z,\pi)}\le C\left(\|u\|_{L^p(Z,\nu)}+\|u\|_{HB^\alpha_{p,q}(Z,\nu)}\right),
    \end{equation*}
    where $C\ge 1$ depends only on $\alpha$, $p$, $q$, $\theta$, $\diam(Z)$, as well as the doubling and codimensional Ahlfors regularity constants.
    
   Furthermore, the bounded linear extension operator $\wtil E:B^{\alpha-\theta/p}_{p,q}(\partial Z,\pi)\to B^\alpha_{p,q}(Z,\nu)$ given by Corollary~\ref{cor:inhom ext bounded domain} is a right-inverse of $T$.  That is, for all $f\in B^{\alpha-\theta/p}_{p,q}(\partial Z,\pi)$, we have that 
    \begin{equation*}
    T\wtil Ef(z)=f(z)
    \end{equation*}
    for $\pi$-a.e.\ $z\in\partial Z$.
    
\end{thm}

To prove Theorem~\ref{thm:Intro gen trace}, we utilize the construction given in \cite{BBS} to realize $(\overline Z,d,\nu)$ as the boundary of a uniform space $(X,d_\eps,\mu_\beta)$, see Section~\ref{sec:hyp fill}, referred to as the uniformized hyperbolic filling of $(\overline Z,d,\nu)$. This space is constructed such that $X$ is uniform in its completion, the measure $\mu_\beta$ is doubling, and $\nu$ is codimensional Ahlfors regular with respect to $\mu_\beta$.  This allows us to define $T$ and $E$ as compositions of operators between the Besov spaces on $\partial Z$, $Z$, and $X$, as given by Theorem~\ref{thm:Trace}, Theorem~\ref{thm:Intro Inhom Trace}, and Theorem~\ref{thm:Extension}.  We use the notation $(Z,d,\nu)$ in the statement of the above theorems to distinguish this space from its hyperbolic filling.  

As mentioned above, Theorem~\ref{thm:Intro gen trace} assumes only that $\nu$ is doubling rather than Ahlfors regular, which is the standing assumption in \cite{Marcos} and \cite{SS}, used primarily to establish the trace theorems there. In \cite[Theorem~2.6]{Marcos}, for example, the trace, or restriction, theorem is obtained via interpolation results for Bessel-type potential spaces in Ahlfors regular metric measure spaces.  As such, it is proven under the assumption that the measure of the space is infinite and that $1<p<\infty$, $1\le q<\infty$, and the regularity parameter $\alpha$ is sufficiently small.  While we assume here that the space has finite measure, Theorem~\ref{thm:Intro gen trace} holds for the full range of regularity parameters $0<\alpha<1-\theta/p$, as well as the cases when $p=1$ or $q=\infty$.  

The trace result \cite[Theorem~1.1]{SS}, which assumes no restriction on the measure of the space, also utilizes the assumption of Ahlfors $Q$-regularity of $\nu$ to ensure that $\pi$-a.e.\ $z\in \partial Z$ is a Lebesgue point of a function $u\in B^\alpha_{p,q}(\overline Z,\nu)$, where $\pi$ is an Ahlfors $d$-regular measure on $\partial Z$, with $0<d<Q$, see \cite[Lemma~3.1]{SS}.  In this manner, the authors show that $Tu=u|_{\partial Z}$ $\pi$-a.e.\ on $\partial Z$ for each $u\in B^\alpha_{p,q}(\overline Z,\nu)$  \cite[Remark~3.3]{SS}.  While many of the estimates used in the proof of their result hold under the weaker assumption of doubling and codimension Ahlfors regularity of $\pi$, it is only shown in this setting that $Tu=u|_{\partial Z}$ for continuous $u\in B^\alpha_{p,q}(\overline Z,\nu)$, as the precise Hausdorff dimension of $\partial Z$ is not known, see \cite[Appendix~A.1]{SS}.  In our result, however, we show that the analogous property \eqref{eq:Intro gen Trace condition} holds for all $u\in B^\alpha_{p,q}(Z,\nu)$, not just continuous functions.    

In Theorem~\ref{thm:Intro gen trace}, the boundedness assumption on $(Z,d,\nu)$ is due to our use of the particular construction of the hyperbolic filling from \cite{BBS}. An alternative hyperbolic filling construction, given via Busemann functions and applicable in the noncompact case, was introduced by Butler in the as-of-yet unpublished papers \cite{Butler1,Butler2}.  In particular, it was shown in \cite{Butler2} that many of the key properties of the uniformized hyperbolic filling from \cite{BBS} hold with this Busemann function construction.  We therefore believe it is likely that Theorem~\ref{thm:Intro gen trace} could be proven without the boundedness assumption on $(Z,d,\nu)$ by using this alternate approach.  However, this investigation is left for future work.

The structure of the paper is as follows: in Section~\ref{section:Prelim}, we discuss the necessary definitions and background notions used throughout the paper.  In Section~\ref{section:Trace}, we prove the trace results Theorem~\ref{thm:Trace} and Theorem~\ref{thm:Intro Inhom Trace}, and Section~\ref{sec:Extension} is devoted to proving the extension results Theorem~\ref{thm:Extension} and Theorem~\ref{thm:InhomExtension}.  Finally, Section~\ref{sec:gen trace thm} is devoted to introducing the necessary background for the construction of the hyperbolic filling and the proof of Theorem \ref{thm:Intro gen trace}.

\section{Preliminaries}\label{section:Prelim}

 In this section, we introduce the necessary definitions and notations used throughout the paper. At times, we let $C>0$ denote a constant which, unless otherwise specified, depends only on the structural constants of the metric measure space, such as the doubling constant, for example.  Its precise value is not of interest to us, and may change with each occurrence, even within the same line.  Furthermore, given quantities $A$ and $B$, we will often use the notation $A\simeq B$ to mean that there exists a constant $C\ge 1$ such that $C^{-1} A\le B\le CA$. Likewise, we use $A\lesssim B$ and $A\gtrsim B$ if the left and right inequalities hold, respectively. Given a metric space $(X,d)$, we will at times use the notation $d(x,A):=\dist(x,A)$ to denote the distance between a point $x\in X$ and a set $A\subset X$.

\subsection{Doubling measures.}

A Borel measure $\mu$ on a metric space $(X,d)$ is said to be \emph{doubling} if there exists a constant $C_\mu\ge 1$ such that 
\[
0<\mu(B(x,2r))\le C_\mu\,\mu(B(x,r))<\infty
\]
for all $x\in X$ and $r>0$.  By iteration, there exists constants $C\ge 1$ and $Q>0$ depending only on $C_\mu$ such that 
\[
\frac{\mu(B(y,r))}{\mu(B(x,R))}\ge C^{-1}\left(\frac{r}{R}\right)^Q
\]
for every $0<r\le R$, $x\in X$, and $y\in B(x,R)$. 

Doubling measures play an important role in analysis on metric measure spaces.  In particular, many classical results, such as Lebesgue differentiation theorem, Hardy-Littlewood maximal function theorem, and the John-Nirenberg inequality can be generalized to the metric setting when the space is equipped with a doubling measure, see for example \cite[Sections 3.4 and 3.5]{HKST} and \cite[Chapter~3]{BB}.  Another important consequence of the doubling property is the following existence of Whitney covers: 

\begin{lem}[{\cite{HKST}, Proposition 4.1.15}]\label{thm:Whitney}
Let $(X,d,\mu)$ be a doubling metric measure space, and let $\Omega\subset X$ be an open set such that $X\setminus\Omega\ne\varnothing$.  Then there exists a countable collection $\mathcal{W}_\Omega=\{B(p_{i,j},r_{i,j})=B_{i,j}\}$ of balls in $\Omega$ so that 
\begin{enumerate}
\item[(i)]$\bigcup_{i,j}B_{i,j}=\Omega$
\item[(ii)] $\sum_{i,j}\chi_{B(p_{i,j},2r_{i,j})}\le 2C_\mu^5$
\item[(iii)] $2^{i-1}<r_{i,j}\le 2^i$ for all $i,$
\item[(iv)] $r_{i,j}=\frac{1}{8}d(p_{i,j}, X\setminus\Omega)$.
\item[(v)] For each $i\in\Z$, $\{p_{i,j}\}_j$ is a $2^i/C$-separated set, and so the collection $\{ B_{i,j}\}_{j\in\Z}$ has bounded overlap.  That is, for each $K\ge 1$, there exists $C_K\ge 1$ so that 
\[
\sum_{j}\chi_{KB_{i,j}}\le C_K.
\]
\end{enumerate}
Furthermore, there exists a Lipschitz partition of unity $\{\phii_{i,j}\}_{i,j}$ subordinate to $\mathcal{W}_\Omega.$  That is, each $\phii_{i,j}$ is $C/r_{i,j}$-Lipschitz continuous with $C$ depending only on $C_\mu$, such that $0\le\phii_{i,j}\le \chi_{2B_{i,j}}$ and $\sum_{i,j}\phii_{i,j}\equiv \chi_{\Omega}$. 
\end{lem}

\begin{defn}\label{def:CDAR}
Let $\theta>0$.  Given an open set $\Omega\subset X$, we say that a measure $\nu$ on $\partial\Omega$ is codimension $\theta$ Ahlfors regular with respect to $\mu$ if there exists a constant $C_\theta\ge 1$ such that for all $z\in\partial\Omega$ and $0<r\le 2\diam(\partial\Omega),$ we have that
\[
\frac{1}{C_\theta}\frac{\mu(B(z,r)\cap\Omega)}{r^\theta}\le\nu(B(z,r)\cap\partial\Omega)\le C_\theta\frac{\mu(B(z,r)\cap\Omega)}{r^\theta}.
\] 
\end{defn}

 For a set $A\subset X$, its codimension $\theta$ Hausdorff measure with respect to $\mu$ is defined as
\[
\Ha^{-\theta}_{\mu}(A):=\lim_{\varepsilon\to 0^+}\, \inf\bigg\lbrace \sum_{i\in I\subset \N}\frac{\mu(B_i)}{\rad(B_i)^\theta}\, 
  :\, A\subset\bigcup_{i\in I}B_i,\, \rad(B_i)\le \varepsilon\bigg\rbrace.
\]
Then we have the following result for codimension $\theta$ Ahlfors regular measures on $\partial\Omega$ (see also \cite[Lemma 2.6]{KLS}).
\begin{lem}\label{lem:codimhausdorff}
    Let $(X,d,\mu )$ be a doubling metric measure space and $\Omega\subset X$ open and such that $\mu|_{\Omega}$ is also doubling. If a measure $\nu$ supported on $\partial \Omega$ is codimension $\theta$ Ahlfors regular with respect to $\mu$ then $\nu\simeq\Ha^{-\theta}_{\mu|_{\overline\Omega}}|_{\partial\Omega}$. In particular, $\Ha^{-\theta}_{\mu|_{\overline\Omega}}|_{\partial\Omega}$ is codimension $\theta$ Ahlfors regular with respect to $\mu$.
\end{lem}
\begin{proof}
Let $\nu$ be codimension $\theta$ Ahlfors regular with respect to $\mu$. Notice that by the doubling condition of $\mu$ and codimensionality we also have $\nu$ doubling (see Remark \ref{rem:doublingboundary} below).

Let $z\in \partial\Omega$, $r>0$. For $\varepsilon >0$ consider a covering of $B(z,r)$ by balls $\{ B_i\}_{i\in I\subset \N}$ in $X$ with $\rad (B_i)<\varepsilon$. Let $I_1\subset I$ the set of indices such that $B_i\cap \partial\Omega\neq\varnothing$, and for each $i\in I$ let $\hat B_i$ with center in $\partial\Omega$ and $\rad (B_i)\leq \rad (\hat B_i)\leq 2\rad (B_i)$, so that $B_i\cap\partial\Omega\subset \hat B_i\subset 3B_i$. Notice that $B(z,r)\cap\partial\Omega \subset\bigcup_{i\in I_1} \hat B_i\cap\partial\Omega$, and then by the doubling condition and codimensionality we have
\begin{align*}
\sum_{i\in I}\frac{\mu|_{\overline\Omega} (B_i)}{\rad (B_i)^\theta}\geq C_d^{-2}\sum_{i\in I_1}\frac{\mu (\hat B_i\cap\Omega )}{\rad (\hat B_i)^\theta}\geq (C_d^2C_\theta)^{-1}\sum_{i\in I_1}\nu (\hat B_i \cap\partial\Omega )\geq (C_d^2C_\theta )^{-1}\nu (B(z,r)\cap\partial\Omega ).
\end{align*}
Taking infimum over all possible coverings and $\varepsilon\to 0$ we get 
$$\nu (B(z,r)\cap\partial\Omega)\lesssim\Ha^{-\theta}_{\mu|_{\overline\Omega}} (B(z,r)\cap\partial\Omega ).$$
Now for each $\eta >0$ there exists $\varepsilon_\eta >0$ such that for all $0<\varepsilon <\varepsilon_\eta$ 
$$\Ha^{-\theta}_{\mu|_{\overline\Omega}} (B(z,r)\cap\partial\Omega )-\eta\leq \inf\bigg\lbrace \sum_{i\in I\subset \N}\frac{\mu|_{\overline\Omega }(B_i)}{\rad(B_i)^\theta}\, 
  :\, B(z,r)\cap\partial\Omega \subset\bigcup_{i\in I}B_i,\, \rad(B_i)\le \varepsilon \bigg\rbrace.$$
Consider also $\varepsilon <r$. By the $5B$-covering lemma, cover $B(z,r)\cap\partial\Omega $ by a collection $\{ B_i\}_{i\in I}$ of balls with radii $\varepsilon$ and centers in $B(z,r)\cap\partial\Omega$ so that $\frac{1}{5}B_i$ are disjoint. Then, since $B_i\subset B(z,2r)$, 
\begin{align*}
\Ha^{-\theta}_{\mu|_{\overline\Omega}} (B(z,r)\cap\partial\Omega )-\eta & \leq 
\sum_{i\in I}\frac{\mu|_{\overline\Omega} (B_i)}{\rad (B_i)^\theta}\leq C_d^{3}\sum_{i\in I}\frac{\mu (\frac{1}{5} B_i\cap\Omega )}{5^\theta\rad (\frac{1}{5} B_i)^\theta}\\
&\leq
\frac{C^3_dC_\theta}{5^\theta}\sum_{i\in I}\nu (\frac{1}{5}B_i\cap\partial\Omega )\leq \frac{C^3_dC_\theta}{5^\theta}\nu (B(z,2r)\cap\partial\Omega).
\end{align*}
Finally, letting $\eta$ go to $0$ we have by the above estimate and the doubling condition on $\nu$
\begin{equation*}\Ha^{-\theta}_{\mu|_{\overline\Omega}} (B(z,r)\cap\partial\Omega )\lesssim \nu (B(z,r)\cap\partial\Omega).\qedhere
\end{equation*}
\end{proof}
\begin{remark}
For a given $\theta >0$ and open set $\Omega$, it follows by definition that all measures on $\partial\Omega$ which are codimension $\theta$ Ahlfors regular with respect to a fixed ambient measure are equivalent. The above lemma shows that $\Ha^{-\theta}_{\mu|_{\overline\Omega}}|_{\partial\Omega}$ is the natural $\theta$-codimensional measure.  Moreover, if this measure is not codimension $\theta$ Ahlfors regular with respect to $\mu$, then there is no such measure for the parameter $\theta$.
\end{remark}

In the work of Jonsson and Wallin on $\R^n$, \cite{JW}, the assumption that $\partial\Omega$ is a $d$-set implies that the $d$-dimensional Hausdorff measure is codimension $(n-d)$ Ahlfors regular with respect to the Lebesgue measure on $\R^n$.  As such, this codimensional assumption on $\nu$ generalizes the $d$-set assumption on $\R^n$.

\begin{remark}\label{rem:doublingboundary}
    If $\mu|_\Omega$ is doubling,  with constant $C_\mu$, and $\nu$ is codimension $\theta$ Ahlfors regular with respect to $\mu$, then $\nu$ is doubling.  Indeed, for $z\in\partial\Omega$ and $0<r\le 2\diam(\partial\Omega)$, we have that
    \[
    \nu(B(z,2r)\cap\partial\Omega)\le C_\theta\frac{\mu(B(z,2r)\cap\Omega)}{(2r)^\theta}\le C_\theta C_\mu\frac{\mu(B(z,r)\cap\Omega)}{(2r)^\theta}\le\frac{C_\mu C_\theta^2}{2^\theta}\nu(B(z,r)\cap\partial\Omega).
    \]
    
\end{remark}

In  Sections~\ref{section:Trace} and \ref{sec:Extension}, we will assume that $(\Omega,d,\mu)$ is a locally compact, non-complete metric measure space, with $\mu$ a doubling measure. The same setting is considered in Section~\ref{sec:gen trace thm}, though with different notation.  Under these assumptions, $\Omega$ is an open set in its completion $\overline\Omega$, and its boundary is given by $\partial\Omega=\overline\Omega\setminus\Omega$.  We extend $\mu$ to $\overline\Omega$ by the zero extension; note that in this setting $\mu$ is doubling on $\overline\Omega$ with the same constant.  Furthermore, we will equip $\partial\Omega$ with a Borel measure $\nu$ which is codimension $\theta$ Ahlfors regular with respect to $\mu$.

\subsection{Uniform domains and Harnack chains}\label{sec:Uniform Domains}

 Throughout this paper, a \emph{curve} is a continuous map $\gamma:[a,b]\to X$. With a slight abuse of notation, we will also use $\gamma$ to denote the trajectory of such a curve.
\begin{defn}\label{def:Uniform domain}
    We say that a domain $\Omega\subset X$ is an  \emph{$A$-uniform domain} for some $A\ge 1$ if for every $x,y\in\Omega$, there exists a curve $\gamma$ in $\Omega$ joining $x$ to $y$ such that 
    \begin{itemize}
        \item[(i)] $\diam(\gamma)\le Ad(x,y)$,  
        \item[(ii)] for each $z\in\gamma$, we have that $\min\{\diam(\gamma_{x,z}),\,\diam(\gamma_{z,y})\}\le Ad(z,X\setminus\Omega)$, where $\gamma_{x,z}$ denotes a fragment of $\gamma$ with endpoints $x$ and $z$ (and analogously for $\gamma_{z,y}$).
    \end{itemize}
\end{defn}

\begin{remark}\label{rem:Uniform Domain}
    In Section~\ref{section:Trace}, we will assume that $\Omega$ is a uniform domain in its completion in order to obtain the desired energy bounds for our trace operator.  In the literature, uniform domains are often defined using the length of curves in place of the diameters used in the above definition, in which case they are sometimes referred to as \emph{length uniform domains}.  Clearly, length uniform domains are uniform domains in the sense of Definition~\ref{def:Uniform domain}. For the study of trace properties of Newton-Sobolev functions, it is often assumed that $\Omega$ is a length uniform domain, see for example \cite{M,Nages-Ryan,BS}, since the standard assumption that the measure is doubling and supports a Poincar\'e inequality ensures that the space is quasiconvex, hence rectifiably connected.

    As our interest is in the nonlocal Besov spaces, we do not assume a Poincar\'e inequality, and so this general setting may include spaces which possess few or possibly zero rectifiable curves.  For this reason, we define uniform domains as above using diameter.  With this definition, the uniform domain property is preserved under a quasisymmetric change in the metric, such as snowflaking for example, and for this reason such definitions are commonly used in the study of Dirichlet forms, see \cite{AChen,Mur,KajMur}.

    To obtain the energy bounds for the Newton-Sobolev trace operators in \cite{M,Nages-Ryan}, the length uniform domain property was used to obtain a suitable chain of balls inside $\Omega$ joining pairs of points on $\partial\Omega$.  As the uniform curves in those settings are rectifiable and thus admit arc-length parametrizations, one can use the Arzel\`a-Ascoli theorem to join each pair of points on $\partial\Omega$ by a uniform curve.  The chain of balls is then obtained by choosing appropriate points along this curve.  We similarly use the uniform domain property to obtain a suitable chain of balls, see Lemma~\ref{lem:uniformcover} below, but as our uniform curves need not be rectifiable, it is not immediate how to extend such curves to the boundary.  For this reason, we obtain our chain of balls in a slightly different manner, using Harnack chains and the following Lemmas~\ref{lem:Unif domain center points} and \ref{lem:Harnack chain} from \cite{AChen}.  
    
\end{remark}

\begin{defn}
    Let $\Omega\subset X$ be a domain and $M\ge 1$.  For $x,y\in\Omega$, an \emph{$M$-Harnack chain} from $x$ to $y$ in $\Omega$ is a sequence of balls $B_1,\dots,B_n$, each contained in $\Omega$, such that $x\in M^{-1}B_{1}$, $y\in M^{-1}B_n$, and $M^{-1}B_i\cap M^{-1}B_{i+1}\ne\varnothing$ for each $i=1,\dots,n-1$. The number $n$ of balls in a Harnack chain is called the \emph{length} of the Harnack chain. 
\end{defn}

\begin{lem}[\cite{AChen}, Lemma 3.4]\label{lem:Unif domain center points} Let $\Omega\subset X$ be an $A$-uniform domain.  For each $z\in\partial\Omega$ and $r>0$ such that $\Omega\setminus B(z,r)\ne\varnothing$, there exists $z_r\in\Omega\cap\partial B(z,r)$ such that $d(z_r,X\setminus\Omega)\ge (2A)^{-1}r$.     
\end{lem}

\begin{lem}[\cite{AChen}, Proposition~2.12]\label{lem:Harnack chain}
Let $(X,d,\mu)$ be a doubling metric measure space.  Let $\Omega\subset X$ be a domain, and for $x,y\in\Omega$, let $\gamma$ be a curve in $\Omega$ joining $x$ to $y$.  Assume that $d(z,X\setminus\Omega)\ge\delta>0$ for all $z\in\gamma$.  Then for any $M>1$ and any $0<r<\delta$, there exists an $M$-Harnack chain $\{B_j:=B(x_j,r)\}_j$ from $x$ to $y$ in $\Omega$, with $x_j\in\gamma$, and with length less than $C(1+Mr^{-1}\diam(\gamma))^\alpha$, where $C$ and $\alpha$ are constants depending only on the doubling constant.    
\end{lem}

\begin{lem}\label{lem:uniformcover}
    Let $(X,d,\mu)$ be a doubling metric measure space, and let $\Omega\subset X$ be an $A$-uniform domain.  Then for each $z,w\in\partial\Omega$, there exists a chain of balls $\{B_k:=B(x_k,r_k)\}_{k\in\Z}$ such that 
    \begin{enumerate}
    \item[(i)]for each $k\in\Z$, we have that $8B_k\subset\Omega$,
    \item[(ii)]$\lim_{k\to\infty}x_k=z$, $\lim_{k\to-\infty}x_k=w$,
    \item[(iii)] for each $k\in\Z$, we have that $2^{-1}B_k\cap 2^{-1}B_{k+1}\ne\varnothing$,
    \item[(iv)] there exists $N\ge 1$ such that for each $k\ge 0$ and each $x\in B_k$, we have that \[
    r_k\simeq d(x,X\setminus\Omega)\simeq d(x,z)\lesssim 2^{-|k|/N}d(z,w),
    \]
    \item[(v)] there exists $N\ge 1$ such that for each $k<0$ and each $x\in B_k$, we have that 
    \[
    r_k\simeq d(x,X\setminus\Omega)\simeq d(x,w)\lesssim 2^{-|k|/N}d(z,w),\]

    \item[(vi)] there exists $C\ge 1$ such that 
    \[
    \sum_{k}\chi_{4B_k}\le C,
    \]
    \end{enumerate}
    The comparison constants, the constant $N$, and the constant $C$ above depend only on $A$ and the doubling constant. 
\end{lem}

\begin{proof}
    Let $z,w\in\partial\Omega$, and let $r_0:=d(z,w)/4$.  By Lemma~\ref{lem:Unif domain center points}, it follows that for each $i\in\N\cup\{0\}
    $, there exists $x_i\in\Omega\cap\partial B(z,r_0/2^i)$ such that
    \begin{equation}\label{eq:xi distance to boundary}
    d(x_i,X\setminus\Omega)\ge\frac{r_0}{(2A)2^i}.
    \end{equation}
    For each $i\in\N\cup\{0\}$, let $\gamma_i$ be an $A$-uniform curve joining $x_i$ to $x_{i+1}$ in $\Omega$. We claim that for all $\zeta\in\gamma_i$, 
    \begin{equation}\label{eq:gammai distance to boundary}
    d(\zeta,X\setminus\Omega)\ge\frac{r_0}{(8A^2)2^i}.
    \end{equation}
    Indeed, from \eqref{eq:xi distance to boundary}, we see that if $\zeta\in\gamma_i\cap B(x_i,r_0/(4A2^i))$, then $d(\zeta,X\setminus\Omega)\ge\frac{r_0}{(4A)2^i}$.  Likewise, if $\zeta\in B(x_{i+1},r_0/(4A2^{i+1}))$, then $d(\zeta,X\setminus\Omega)\ge\frac{r_0}{(8A)2^i}$.  For the remaining case when 
    \[    \zeta\in\gamma_i\setminus\left(B\left(x_i,\frac{r_0}{(4A)2^i}\right)\cup B\left(x_{i+1},\frac{r_0}{(4A)2^{i+1}}\right)\right),
    \]
    it follows from the $A$-uniform property of $\gamma_i$ that 
    \begin{align*}
        d(\zeta,X\setminus\Omega)\ge A^{-1}\min\{\diam(\gamma_{i,x_i,\zeta}),\,\diam(\gamma_{i,\zeta,x_{i+1}})\}\ge\frac{r_0}{(8A^2)2^i}.
    \end{align*}
    Thus, \eqref{eq:gammai distance to boundary} holds in each case, proving the claim.

    Letting 
    \[
    r_i:=\frac{r_0}{(64A^2)2^i},
    \]
    we then apply Lemma~\ref{lem:Harnack chain} to obtain a $2$-Harnack chain $\{B_{i,j}:=B(x_{i,j},r_i)\}_{j=1}^{n_i}$ from $x_i$ to $x_{i+1}$ in $\Omega$, with $x_{i,j}\in\gamma_i$, and such that
    \begin{equation}\label{eq:Harnack chain length}
        n_i\le C(1+2r_i^{-1}\diam(\gamma_i))^\alpha\le C(1+192 A^3)^\alpha=:N.
    \end{equation}
    Here the last inequality follows from the coarse estimate that $\diam(\gamma_i)\le(3A/2)r_0/2^i$.  We claim that for each $j=1,\dots,n_i$, the following holds for all $x\in B_{i,j}$, with comparison constants depending only on $A$ and the doubling constant:
    \begin{align}\label{eq:Bij distance to boundary}
        r_i\simeq d(x,X\setminus\Omega)\simeq d(x,z)\simeq 2^{-i}d(z,w).
    \end{align}
    Indeed, $r_i\simeq 2^{-i}d(z,w)$ by the definition of $r_i$ and $r_0$.  Moreover, we have that 
    \begin{align}\label{eq:ri comparison}
        d(x,X\setminus\Omega)\le d(x,z)\le d(x,x_{i,j})+d(x_{i,j},x_{i+1})+d(x_{i+1},z)\le r_i+\diam(\gamma_i)+r_0/2^{i+1}\lesssim r_i,
    \end{align}
    and by \eqref{eq:gammai distance to boundary}, we have that
    \[
    d(x,X\setminus\Omega)\ge d(x_{i,j},X\setminus\Omega)-d(x_{i,j},x)\gtrsim r_i.
    \]
    Thus, \eqref{eq:Bij distance to boundary} holds, proving the claim.

    By relabeling the collection $\bigcup_{i=0}^\infty\{B_{i,j}\}_{j=1}^{n_i}=\{B_k:=B(x_k,r_k)\}_{k\ge 0}$ with the appropriate ordering, it follows from the construction, the definition of $2$-Harnack chains, and \eqref{eq:gammai distance to boundary}, that claims $(i)$, $(ii)$, and $(iii)$ of the lemma hold (when $k\ge 0$).  Furthermore, if $B_k=B(x_{i,j},r_i)$ for some $i,j$, then from \eqref{eq:Bij distance to boundary} and \eqref{eq:Harnack chain length}, it follows that $k\le Ni$.  Thus, $-i\le -k/N$, and so claim $(iv)$ holds.
    
    To define the collection of balls $\{B_k\}_{k<0}$, we repeat the argument above, this time using Lemma~\ref{lem:Unif domain center points} to obtain $y_i\in\Omega\cap\partial B(w,r_0/2^i)$ such that $d(y,X\setminus\Omega)\ge r_0/(2A2^i)$ for each $i\in\N\cup\{0\}$.  We then use uniform curves joining $y_i$ to $y_{i+1}$ and Lemma~\ref{lem:Harnack chain} to obtain the analogous collection of balls $\{B(y_{i,j},r_i)\}_{j=1}^{n_i}$ via $2$-Harnack chains, with length bounded by $N$, along the uniform curves.  We also join $x_0$ to $y_0$ by an $A$-uniform curve, and then apply Lemma~\ref{lem:Harnack chain} to likewise obtain a $2$-Harnack chain $\{B(z_{0,j},r_0/(64A^2))\}_j$ from $x_0$ to $y_0$ of length bounded by $N$.  We then relabel the collection
    \[
    \{B(z_{0,j},r_0/(64A^2))\}_j\cup\bigcup_{i=0}^\infty\{B(y_{i,j},r_i)\}_{j=1}^{n_i}=\{B_k:=B(x_k,r_k)\}_{k<0},
    \]
    and by the above arguments, it follows that the collection $\{B_k\}_{k\in\Z}$ satisfies claims $(i)$-$(v)$ of the Lemma.

    It remains to prove claim $(vi)$.  To this end, we denote 
    \[
    \B_i^+:=\{B(x_{i,j},r_i)\}_{j},\quad\B_i^-:=\{B(y_{i,j},r_i)\}_j,\quad\B_0:=\{B(z_{i,j},r_0/(64A^2))\}_j.
    \]
    Then $\{B_k\}_{k\in\Z}$ is just a relabeling of $\B_0\cup\bigcup_{i=0}^\infty(\B_i^+\cup\B_i^-)$.  Let $B:=B(x_{i,j},r_i)$ be a ball in the collection $\B_i^+$, and let $x\in 4B$.  Then, from \eqref{eq:gammai distance to boundary}, it follows that 
    \begin{align*}
    d(x,X\setminus\Omega)\ge d(x_{i,j},X\setminus\Omega)-d(x,x_{i,j})\ge\frac{r_0}{(8A^2)2^i}-\frac{r_0}{(16A^2)2^i}=\frac{r_0}{(16A^2)2^i}.
    \end{align*}
    Likewise, similar to \eqref{eq:ri comparison}, we have that
    \begin{align*}
        d(x,X\setminus\Omega)\le d(x,z)\le\frac{r_0}{(16A^2)2^i}+(3A/2)\frac{r_0}{2^i}+\frac{r_0}{2^{i+1}}=\left((16A^2)^{-1}+3A/2+1/2\right)\frac{r_0}{2^i}.
    \end{align*}
    This shows that there exists $C\ge1$ depending only on $A$ such that if $i,i'\in\N\cup\{0\}$ with $|i-i'|\ge C$, then
    \[
    4B\cap 4B'=\varnothing
    \]
    for all $B\in\B^+_i$ and $B'\in \B^+_{i'}$, and by symmetry for all $B\in\B^-_{i}$ and $B'\in\B^-_{i'}$. Similar arguments show that there exists (a possibly different) $C\ge 1$, depending only on $A$, such that if $i,i'\in\N\cup\{0\}$ are such that $i,i'\ge C$, then for all $B\in\B_i^+$, $B'\in\B^-_{i'}$, $B''\in\B_0$, we have
    \[
    4B\cap 4B'=4B\cap 4B''=4B'\cap4B''=\varnothing.
    \]
    Thus, for a given $B_k$, if
    \[
    4B_{k'}\cap 4B_k\ne\varnothing,
    \]
    then there are at most $C$ of the collections $\B_0$, $\B_i^{+}$, $\B_{i}^-$, $i=0,1,2,\dots$ in which $B_{k'}$ may belong.  As there are at most $N$ balls in each of these collections, with $N$ depending only on $A$ and the doubling constant, it follows that claim $(vi)$ holds. 
\end{proof}

\subsection{Besov spaces.}\label{section:Besov} 
Let $1\leq p,q<\infty$ and $0<\alpha<1$. Given a function $f\in L^1_\loc(X,\mu)$, we define its Besov  energy by
\[
\Vert f\Vert_{HB^\alpha_{p,q}(X,\mu )}:= \left(\int_0^\infty\left(\int_X\fint_{B(x,t)}|f(y)-f(x)|^pd\mu(y)d\mu(x)\right)^{q/p}\frac{dt}{t^{\alpha q+1}}\right)^{1/q}.\]
As usual, we define the homogeneous Besov space $HB^\alpha_{p,q}(X,\mu)$ as the class of functions in $L^1_\loc(X,\mu)$ for which this energy is finite, while the inhomogeneous Besov space is defined by $B^\alpha_{p,q}(X,\mu ):=HB^\alpha_{p,q}(X,\mu )\cap L^p(X,\mu )$, and is equipped with the norm 
\[
\|u\|_{B^\alpha_{p,q}(X,\mu)}:=\|u\|_{HB^\alpha_{p,q}(X,\mu)}+\|u\|_{L^p(X,\mu)}.
\]
We follow the ideas in the proof of \cite{GKS} in order to obtain the following equivalent seminorm for the Besov energy.

\begin{lem}
Let $(X,d,\mu )$ be a metric measure space with $\mu$ a doubling measure. 
	Let $1\leq p,q<\infty$ and $0<\alpha<1$. Given a function $f\in L^1_\loc(X,\mu)$, then
\begin{equation}\label{eq:SumForm}
	\|f\|^q_{HB^\alpha_{p,q}( X ,\mu)}\simeq \sum_{k\in\Z}2^{-k\alpha q}\left(\int_ X \fint_{B(x,2^k)}|f(y)-f(x)|^pd\mu(y)d\mu(x)\right)^{q/p},
\end{equation}
where the comparison constants depend only on $\alpha$, $p$, $q$ and $C_\mu$. Moreover, if $f\in L^p( X,\mu )$ then
\begin{equation}\label{eq:SumFormLp}
	\|f\|^q_{HB^\alpha_{p,q}( X ,\mu)}\lesssim \| f\|^q_{L^p( X ,\mu )}+\sum_{k\in\N}2^{k\alpha q}\left(\int_ X \fint_{B(x,2^{-k})}|f(y)-f(x)|^pd\mu(y)d\mu(x)\right)^{q/p},
\end{equation}
where the comparison constants depend only on $\alpha$, $p$, $q$ and $C_\mu$.
\end{lem}
\begin{proof}
	By the doubling property we obtain \eqref{eq:SumForm}  following the steps in \cite[Theorem 5.2]{GKS}.  To prove \eqref{eq:SumFormLp}, consider $k\in \N$.  Then
$$
	\int_ X\fint_{B(x,2^k)}|f(x)-f(y)|^pd\mu (y)d\mu (x) \lesssim 
	\int_ X |f(x)|^pd\mu (x)+\int_ X \fint_{B(x,2^k)}|f(y)|^pd\mu (y)d\mu(x).
$$

By Tonelli's theorem and doubling, we have that
\begin{align*}
    \int_X\fint_{B(x,2^k)}|f(y)|^pd\mu(y)d\mu(x)&=\int_X\int_{X}|f(y)|^p\frac{\chi_{B(x,2^k)}(y)}{\mu(B(x,2^k))}d\mu(y)d\mu(x)\\
    &=\int_X|f(y)|^p\int_{X}\frac{\chi_{B(y,2^k)}(x)}{\mu(B(x,2^k))}d\mu(x)d\mu(y)\lesssim\int_X|f(y)|^pd\mu(y).
\end{align*}
Since $\sum_{k\in\N}2^{-k\alpha q}$ converges, \eqref{eq:SumFormLp} holds by splitting the sum in \eqref{eq:SumForm} into $k>0$ and $k\le 0$.
\end{proof}

We now record the following sum-rearrangement lemma from \cite{HIT}, which we will use frequently in establishing the energy bounds for our trace and extension operators. 

\begin{lem}[\cite{HIT}, Lemma 3.1]\label{lem:Sums}
    Let $1<a<\infty$, $0<b<\infty$, and $c_i\ge 0$, $i\in\Z$.  Then there is a constant $C=C(a,b)$ such that 
    \[
    \sum_{i\in\Z}\left(\sum_{j\in\Z}a^{-|j-i|}c_j\right)^b\le C\sum_{j\in \Z}c_j^b.
    \]
\end{lem}

To close this section, we provide the definition of the Besov class for $q=\infty$, which will be treated separately within the proofs of the main results.

\begin{align*}
\Vert f\Vert_{HB^\alpha_{p,\infty}(X,\mu )}&:= \sup_{t>0}\frac{1}{t^{\alpha}}\left(\int_X\fint_{B(x,t)  }|f(y)-f(x)|^pd\mu(y)d\mu(x)\right)^{1/p}\\
&\simeq \sup_{k\in\Z}2^{-k\alpha }\left(\int_X \fint_{B(x,2^k)  }|f(y)-f(x)|^pd\mu(y)d\mu(x)\right)^{1/p},
\end{align*}
Moreover, if $f\in L^\infty (X ,\mu )$ then
\[
\Vert f\Vert_{HB^\alpha_{p,\infty}(X,\mu )}\lesssim \Vert f\Vert_{L^\infty (X ,\mu )}+\sup_{k\in\N}2^{k\alpha }\left(\int_X \fint_{B(x,2^{-k})  }|f(y)-f(x)|^pd\mu(y)d\mu(x)\right)^{1/p},
\]
The seminorm equivalences above are obtained analogously to those in \eqref{eq:SumForm} and \eqref{eq:SumFormLp}.

\section{Proof of Trace Theorems~\ref{thm:Trace} and \ref{thm:Intro Inhom Trace}}\label{section:Trace}

In this section, we assume that $(\Omega,d,\mu)$ is a locally compact, non-complete metric measure space, with $\mu$ a doubling measure, such that $\Omega$ is a uniform domain in its completion $\overline\Omega$, see Definition~\ref{def:Uniform domain}. We extend $\mu$ to $\overline\Omega$ by the zero extension, and note that the metric measure space $(\overline\Omega,d,\mu)$ is still doubling with the same constant. Due to this zero extension, it follows that $HB^\alpha_{p,q}(\Omega,\mu)=HB^\alpha_{p,q}(\overline\Omega,\mu)$; by a slight abuse of notation, we will at times use the spaces interchangeably.  We equip  $\partial\Omega:=\overline\Omega\setminus\Omega$, the boundary of $\Omega$, with a Borel measure $\nu$ which is codimension $\theta$ Ahlfors regular with respect to $\mu$ for some $\theta>0$. Recall that this implies that $\nu$ is doubling, see Remark~\ref{rem:doublingboundary}.  The arguments in this section are inspired by those of \cite{M,Nages-Ryan}, where analogous trace results were proven for the Newton-Sobolev class.

\subsection{Existence of traces and $L^p$ bounds.} 
Fix $R>0$ and let $\theta>0$ be as above. We define the following restricted fractional maximal function following \cite{M},
\[
M^R_\theta f(z)=\sup_{0<r<R}r^\theta\fint_{B(z,r)}|f|d\mu,
\]
which maps functions in $L^1_{\loc}(\Omega)$ to the space of lower semicontinuous functions on $\partial\Omega$.  The following lemma from \cite{M} was proven in greater generality; for the reader's convenience, we include the proof of the case relevant to us, that $M^R_\theta$ is bounded from $L^1(\Omega)$ to weak-$L^1(\partial\Omega)$:

\begin{lem}[\cite{M}, Lemma 4.2]\label{lem:Weak1,1}
There exists a constant $C$, depending only on $\theta$, $C_\mu$ and $C_\theta$, such that for all $\lambda>0$ and $f\in L^1(\Omega,\mu)$,
\[
\nu(\{z\in\partial\Omega:M^R_\theta f(z)>\lambda\})\le \frac{C}{\lambda}\int_\Omega|f|d\mu.
\]
\end{lem}

\begin{proof}
    Let $E_\lambda:=\{z\in\partial\Omega:M^R_\theta f(z)>\lambda\}$.  Then for each $z\in E_\lambda$, there exists $B_z:=B(z,r_z)$ such that $r_z^\theta\fint_{B_z}|f|d\mu>\lambda$, and so we have that 
    \[
    \frac{\mu(B_z)}{r_z^\theta}\le\frac{1}{\lambda}\int_{B_z}|f|d\mu.
    \]
    Since $r_z<R$ for all $z\in E_\lambda$, we use the 5-covering lemma to obtain a disjoint countable subcollection $\{B_i:=B(z_i,r_{z_i})\}_i$ such that 
    \[
    \bigcup_{z\in E_\lambda} B_z\subset\bigcup_i5B_i.
    \]
    Therefore, it follows from doubling, codimensionality of $\nu$, and the disjointness of $\{B_i\}_i$ that 
    \begin{align*}
        \nu(E_\lambda)\le\nu\left(\bigcup_i 5B_i\cap\partial\Omega\right)&\le\sum_i\nu(5B_i\cap\partial\Omega)\\
        &\lesssim\sum_i\frac{\mu(B_i)}{r_{z_i}^\theta}\le\frac{1}{\lambda}\sum_i\int_{B_i}|f|d\mu\le\frac{1}{\lambda}\int_{\Omega}|f|d\mu.\qedhere
    \end{align*}
   
\end{proof}

Before studying the trace operator we need the following lemma that will ensure the trace is well defined. 
\begin{lem}\label{lem:L^1 loc on boundary}
Let $1\le p<\infty$, $1\leq q\leq\infty$, $0<\alpha<1$, and let $u\in HB^\alpha_{p,q}(\Omega,\mu)$.  Then for each $z\in\partial\Omega$ and $r>0$, we have that 
\[
\int_{B(z,r)}|u|d\mu<\infty.
\]    
\end{lem}
\begin{proof}
    Since $\|u\|_{HB^\alpha_{p,q}(\Omega,\mu)}<\infty$, it follows that 
    \begin{align*}
    \infty>\int_{B(z,r) }\fint_{B(x,t) }|u(x)-u(y)|^pd\mu(y)d\mu(x)
    \end{align*}
    for $\Leb$-a.e.\ $0<t<\infty$.  Since $u\in L^1_\loc(\Omega,\mu)$, we then have that $|u(x)|<\infty$ and
    \[
    \infty>\fint_{B(x,t) }|u(x)-u(y)|^pd\mu(y)
    \]
 for $\mu$-a.e.\ $x\in B(z,r) $.  Choosing such $t\ge 2r$ and $x\in B(z,r) $, the conclusion follows, as 
    \begin{align*}
        \infty>\left(\fint_{B(x,t) }|u(x)-u(y)|^pd\mu(y)\right)^{1/p}&\ge\fint_{B(x,t) }|u(x)-u(y)|d\mu(y)\\
        &\ge\frac{1}{\mu(B(x,t) )}\left(\int_{B(z,r) }|u(y)|d\mu(y)-|u(x)|\right).\qedhere
    \end{align*}
\end{proof}
    
We note that the assumption that $\Omega$ is uniform in its completion is not needed to prove the following result.  However, this assumption will be used in the following subsection to attain the appropriate energy bounds.

\begin{prop}\label{prop:tracelimit}
Let $1\le p<\infty$, $1\leq q\leq\infty $ and $0<\alpha<1$ be such that $\alpha-\theta/p>0$.  Then there exists a linear operator $T: HB^\alpha_{p,q}(\Omega,\mu)\rightarrow L^p(\partial\Omega,\nu)$ such that 
\begin{equation}\label{eq:TraceCondition}
\lim_{r\to 0^+}\fint_{B(z,r) }|u-Tu(z)|^pd\mu=0\text{ for }\nu\text{-a.e\ }z\in\partial\Omega.
\end{equation}
If in addition $\partial\Omega$ is bounded, then 
\begin{equation}\label{eq:Trace Lp bound}
\|Tu\|_{L^p(\partial\Omega,\nu)}\lesssim \diam(\partial\Omega)^{-\theta/p}\|u\|_{L^p(\Omega,\mu)}+\diam(\partial\Omega)^{\alpha -\theta/p}\|u\|_{HB^{\alpha}_{p,q}(\Omega,\mu)}.
\end{equation}
\end{prop}

\begin{proof}
We first consider the case $1\le q<\infty$.  The case $q=\infty$ holds with natural modifications briefly described at the end of the proof.  For each $r>0$, choose a maximal $r$-separated subset $\{z_i^r\}_{i\in I_r}$ of $\partial\Omega$.  For each $i\in I_r$, set $B_i^r:=B(z_i^r,r)$, and let $U_i^r:=B(z_i^r,r)\cap\partial\Omega$.  Since $\mu$ is doubling, there exists a Lipschitz partition of unity $\{\phii_i^r\}_{i\in I_r}$ subordinated to $\{B_i^r\}_{i\in I_r}$.  That is, each $\phii_i^r$ is $C/r$-Lipschitz, $0\le \phii_i^r\le\chi_{2B_i^r}$, and for all $x\in\bigcup_{i\in I_r} B_i^r$, we have that $\sum_{i\in I_r}\phii_i^r(x)=1$.  For proof of these facts, see \cite[Lemmas B.7.3 and B.7.4]{Gromov}, for example. 

For each $u\in HB_{p,q}^\alpha(\Omega,\mu)$, we then set 
\[
u_{B^r_i}:=\fint_{B_i^r }u\,d\mu,
\]
and define 
\[
T_ru:=\sum_{i\in I_r}u_{B_i^r}\phii_i^r|_{\partial\Omega}.
\]
Note that $T_ru$ well-defined by Lemma~\ref{lem:L^1 loc on boundary} and locally Lipschitz, hence $\nu$-measurable.  We claim that for any sequence $r_k\to 0^+$, the sequence $\{T_{r_k}u\}_k$ is Cauchy in $L^p(\partial\Omega,\nu)$.  Toward this, consider $r,R>0$ such that $0<R/2<r\le R$.  Since $\mu$ is doubling, we have
\begin{align}\label{eq:intersected comparable measure}
\mu(B(z,r))\simeq \mu(B(z,R)).
\end{align}
Setting
\[
\beta:=\frac{\alpha p+\theta}{2},\quad\sigma:=\alpha-\beta/p,
\]
we note that $\sigma>0$ since $\alpha>\theta/p$.  Using the properties of the partition of unity, we then have that 

\begin{align*}
    \|T_Ru-T_ru\|^p_{L^p(\partial\Omega,\nu)}&\le\sum_{i\in I_r}\int_{U_i^r}\left|\sum_{j\in I_r}u_{B_j^r}\phii_j^r(z)-\sum_{k\in I_R}u_{B_k^R}\phii_k^R(z)\right|^pd\nu(z)\\
    &\le\sum_{i\in I_r}\int_{U_i^r}\left(\sum_{j\in I_r}|u_{B_j^r}-u_{B_i^r}|\phii_j^r(z)+\sum_{k\in I_R}|u_{B_i^r}-u_{B_k^R}|\phii_k^R(z)\right)^pd\nu(z)\\
    &\lesssim\sum_{i\in I_r}\int_{U_i^r}\fint_{B(z,3R)}\fint_{B(z,3R)}|u(x)-u(y)|^pd\mu(y)d\mu(x)d\nu(z)\\
    &\lesssim\int_{\partial\Omega}\fint_{B(z,3R)}\fint_{B(z,3R)}|u(x)-u(y)|^pd\mu(y)d\mu(x)d\nu(z).
\end{align*}
To obtain the last two inequalities, we have used the fact that if $z\in U_i^r$, and $\phii_k^R(z)\ne 0$, then $B_k^R\subset B(z,3R)$ (and similarly if $\phii_j^r(z)\ne 0$), along with \eqref{eq:intersected comparable measure} and bounded overlap of the collection of balls, due to the doubling property.

Furthermore, we have
\begin{align*}
    \|T_Ru-&T_ru\|_{L^p(\partial\Omega,\nu)}^p\\
    &\lesssim R^{\beta}\int_{\partial\Omega}\fint_{B(z,3R) }\int_{B(z,3R) }\frac{|u(x)-u(y)|^p}{d(x,y)^{\beta}\mu(B(x,d(x,y)) )}d\mu(y)d\mu(x)d\nu(z)\\
    &\lesssim R^{\beta-\theta}\int_{\partial\Omega}\int_{B(z,3R) }\int_{B(z,3R) }\frac{|u(x)-u(y)|^p\,d\mu(y)d\mu(x)}{d(x,y)^{\beta}\mu(B(x,d(x,y)))\nu(B(z,R)\cap\partial\Omega)}d\nu(z).
\end{align*}

Letting $\Omega_{3R}:=\{x\in\Omega:d(x,\partial\Omega )<3R\}$, we have by Tonelli's theorem that 
\begin{align*}
\|T_Ru&-T_ru\|_{L^p(\partial\Omega,\nu)}^p\nonumber\\
    &\lesssim R^{\beta-\theta}\int_{\Omega_{3R}}\int_{B(x,6R) }\int_{B(x,3R)\cap\partial\Omega}\frac{|u(x)-u(y)|^pd\nu(z)d\mu(y)d\mu(x)}{d(x,y)^{\beta}\mu(B(x,d(x,y)) )\nu(B(z,R)\cap\partial\Omega)}\nonumber\\
    &\lesssim R^{\beta-\theta}\int_{\Omega}\int_{B(x,6R) }\frac{|u(x)-u(y)|^p}{d(x,y)^{\beta}\mu(B(x,d(x,y)) )}d\mu(y)d\mu(x).
\end{align*}
Let $N_R\in\Z$ be such that $2^{N_R-1}\le 6R<2^{N_R}$. For each $x\in \Omega$ and $m\in\Z$, consider the annulus $A_m(x):= B(x,2^{-m})\backslash B(x,2^{-m-1})$. We then have that 
\begin{align*}
\|T_Ru-T_ru\|_{L^p(\partial\Omega,\nu)}\nonumber
    &\lesssim \left( R^{\beta-\theta}\int_{\Omega}\sum^\infty_{m=-N_R}\int_{A_m(x)}\frac{|u(x)-u(y)|^p}{d(x,y)^{\beta}\mu(B(x,d(x,y)) )}d\mu(y)d\mu(x)\right)^{1/p}\nonumber\\
    &\lesssim\left( R^{\beta-\theta}\sum^\infty_{m=-N_R}2^{m\beta}\int_{\Omega}\int_{A_m(x)}\frac{|u(x)-u(y)|^p}{\mu(B(x,2^{-m}) )}d\mu(y)d\mu(x)\right)^{1/p}\nonumber\\
    &\le\left( R^{\beta-\theta}\sum^\infty_{m=-N_R}2^{m\beta}\int_{\Omega}\fint_{B(x,2^{-m}) }|u(x)-u(y)|^pd\mu(y)d\mu(x)\right)^{1/p}\nonumber\\
    &\le R^{(\beta-\theta)/p}\ \sum^\infty_{m=-N_R}2^{m\beta/p}\left(\int_{\Omega}\fint_{B(x,2^{-m}) }|u(x)-u(y)|^pd\mu(y)d\mu(x)\right)^{1/p}\nonumber.
\end{align*}
Here in the last inequality, we have used the fact that for any $a_i>0$, $i\in\N$, and $0<b\leq 1$, it follows that
\[\left(\sum_{i\in\N}a_i\right)^b\leq \sum_{i\in\N}a_i^b.
\]
By Hölder's inequality, and our choices of $\beta$, $\sigma$, and $N_R$, we then obtain
\begin{align*}
    \|T&_Ru-T_ru\|_{L^p(\partial\Omega,\nu)}\\
    &\lesssim R^{\frac{\beta-\theta}{p}}\left(\sum_{m=-N_R}^\infty 2^{\frac{-m\sigma q}{q-1}}\right)^{\frac{q-1}{q}}\left(\sum^\infty_{m=-N_R}2^{m\left(\frac{\beta}{p}+\sigma\right)q}\left(\int_{\Omega}\fint_{B(x,2^{-m}) }|u(x)-u(y)|^pd\mu(y)d\mu(x)\right)^{\frac{q}{p}}\right)^{\frac{1}{q}}\\
    &\simeq R^{(\beta-\theta)/p+\sigma}\left(\sum^\infty_{m=-N_R}2^{m\alpha q}\left(\int_{\Omega}\fint_{B(x,2^{-m}) }|u(x)-u(y)|^pd\mu(y)d\mu(x)\right)^{q/p}\right)^{1/q}\\
    &\le R^{\alpha-\theta/p}\|u\|_{HB^\alpha_{p,q}(\Omega,\mu)}.
\end{align*}
Now, for any $0<r\le R$, there exists $N\in\N$ such that $2^{-N}R<r\le 2^{-N+1}R$.  Thus, we have that 
\begin{align}\label{eqn:L^pCauchy}
\|T_Ru-T_ru\|_{L^p(\partial\Omega,\nu)}&\le\sum_{k=1}^N\|T_{2^{-k+1}R}u-T_{2^{-k}R}u\|_{L^p(\partial\Omega,\nu)}\nonumber\\
    &\lesssim\|u\|_{HB^\alpha_{p,q}(\Omega,\mu)}\sum_{k=1}^N(2^{-k}R)^{\alpha-\theta/p}\lesssim R^{\alpha-\theta/p}\|u\|_{HB^\alpha_{p,q}(\Omega,\mu)}\to 0
\end{align}
as $R\to 0^+$, since $u\in HB^\alpha_{p,q}(\Omega,\mu)$.  Hence, for any sequence $r_k\to 0^+$, $\{T_{r_k}u\}_k$ is Cauchy in $L^p(\partial\Omega,\nu)$, and so there exists $Tu\in L^p(\partial\Omega,\nu)$ such that $T_ru\to Tu$ in $L^p(\partial\Omega,\nu)$ as $r\to 0^+$.  Passing to a subsequence if necessary, we have that $T_ru(z)\to Tu(z)$ for $\nu$-a.e.\ $z\in\partial\Omega$.   

Now, fix $0<R\leq 1$ and define
\[
F(x):=\int_{B(x,R) }\frac{|u(x)-u(y)|^p}{d(x,y)^{\beta}\mu(B(x,d(x,y)))}d\mu(y).
\]
Then, using a similar argument via H\"older's inequality and doubling as above, we have that 
\begin{align*}
    \left(\int_\Omega |F|d\mu\right)^{1/p}&
    = \left(\int_\Omega \int_{B(x,R) }\frac{|u(x)-u(y)|^p}{d(x,y)^{\beta}\mu(B(x,d(x,y)))}d\mu(y)d\mu(x)\right)^{1/p}\\
    &\le\left(\sum_{m\geq 0}\int_\Omega\int_{B(x,2^{-m+1})\setminus B(x,2^{-m}) }\frac{|u(x)-u(y)|^p}{d(x,y)^{\beta}\mu(B(x,d(x,y)))}d\mu(y)d\mu(x)\right)^{1/p}\\
    &\lesssim\left(\sum_{m\geq 0}2^{m\beta}\int_\Omega\fint_{B(x,2^{-m}) }|u(x)-u(y)|^pd\mu(y)d\mu(x)\right)^{1/p}\\
   &\le\sum_{m\geq 0}2^{m\beta/p}\left(\int_\Omega\fint_{B(x,2^{-m}) }|u(x)-u(y)|^pd\mu(y)d\mu(x)\right)^{1/p}\\
    &\lesssim \left(\sum_{m\geq 0}2^{m(\beta/p+\sigma)q}\left(\int_\Omega\fint_{B(x,2^{-m}) }|u(x)-u(y)|^pd\mu (y)d\mu (x)\right)^{q/p}\right)^{1/q}\\
    &\lesssim\| u\|_{HB^\alpha_{p,q}(\Omega ,\mu)}.
\end{align*}
Since $u\in HB^\alpha_{p,q}(\Omega,\mu)$, it follows that $F\in L^1(\Omega,\mu)$, and so by Lemma~\ref{lem:Weak1,1}, we have that $M^R_\theta F(z)<\infty$ for $\nu$-a.e.\ $z\in\partial\Omega$. 

Let $z\in\partial\Omega$ such that $M^R_{\theta}F(z)<\infty$ and $T_ru(z)\to Tu(z)$ as $r\to 0^+$.  This holds for $\nu$-a.e.\ $z\in\partial\Omega$.  Let $0<r<R$.  Then, by properties of the partition of unity, \eqref{eq:intersected comparable measure}, our definition of $\beta$, and since $\alpha>\theta/p$, we have that 
\begin{align*}
\fint_{B(z,r) }|u-&T_ru(z)|^pd\mu=\fint_{B(z,r) }\left|u(x)-\sum_{i\in I_r}u_{B_i^r}\phii_i^r(z)\right|^pd\mu(x)\\   &\lesssim\fint_{B(z,r) }\fint_{B(z,3r) }|u(x)-u(y)|^pd\mu(y)d\mu(x)\\
    &\lesssim r^{\beta}\fint_{B(z,r) }\int_{B(z,3r) }\frac{|u(x)-u(y)|^p}{d(x,y)^{\beta}\mu(B(x,d(x,y)) )}d\mu(y)d\mu(x)\\
    &\le r^{\beta-\theta}r^\theta\fint_{B(z,r) }\int_{B(z,R) }\frac{|u(x)-u(y)|^p}{d(x,y)^{\beta}\mu(B(x,d(x,y)) )}d\mu(y)d\mu(x)\\
    &= r^{(\alpha p-\theta)/2}\left(r^\theta\fint_{B(z,r) }|F(x)|d\mu(x)\right)\le r^{(\alpha p-\theta)/2}M^R_{\theta}F(z)\to 0,
\end{align*}
as $r\to 0^+$.  Hence, we have that 
\begin{align*}
\fint_{B(z,r) }|u-Tu(z)|^pd\mu\lesssim\fint_{B(z,r) }|u-T_ru(z)|^pd\mu+|T_ru(z)-Tu(z)|^p\to 0
\end{align*}
as $r\to 0^+$, which gives us \eqref{eq:TraceCondition}.

Assume now that $\partial\Omega$ is bounded, and let $R=2\diam (\partial\Omega)$. Then, 
\begin{align*}
\|Tu\|_{L^p(\partial\Omega,\nu)}&\le\lim_{r\to 0^+}\|T_ru-T_Ru\|_{L^p(\partial\Omega,\nu)}+\|T_Ru\|_{L^p(\partial\Omega,\nu)},
\end{align*}
By \eqref{eqn:L^pCauchy}, we have that
\[
\|T_ru-T_Ru\|_{L^p(\partial\Omega,\nu)}\lesssim \diam(\partial\Omega)^{\alpha -\theta/p}\|u\|_{HB^\alpha_{p,q}(\Omega,\mu)}.
\]
By properties of the partition of unity, as well as bounded overlap of $\{B_i^R\}_{i\in I_R}$, we also have that 
\begin{align*}
\|T_Ru\|^p_{L^p(\partial\Omega,\nu)}&\lesssim\int_{\partial\Omega}\fint_{B(z,3R) }|u|^pd\mu d\nu(z)\\
&\le\|u\|^p_{L^p(\Omega,\mu)}\int_{\partial\Omega}\frac{1}{\mu(B(z,3R) )}d\nu(z)\\
    &\lesssim\|u\|^p_{L^p(\Omega,\mu)}\int_{\partial\Omega}\frac{1}{R^\theta\nu(B(z,R)\cap\partial\Omega)}d\nu(z)=(2\diam(\partial\Omega))^{-\theta}\|u\|^p_{L^p(\Omega,\mu)}.
\end{align*}
Hence, it follows that 
\[
\|Tu\|_{L^p(\partial\Omega,\nu)}\lesssim \diam(\partial\Omega)^{-\theta/p}\|u\|_{L^p(\Omega,\mu)}+\diam(\partial\Omega)^{\alpha -\theta/p}\|u\|_{HB^{\alpha}_{p,q}(\Omega,\mu)},
\]
which gives us \eqref{eq:Trace Lp bound}.

For the case $q=\infty $ we notice that
\begin{align*}
\sum_{m=-N_R}^\infty 2^{m\beta/p}\left(\int_\Omega\fint_{B(x,2^{-m}) }|u(x)-u(y)|^pd\mu (y)d\mu (x)\right)^{1/p} &\lesssim 2^{-N_R(\beta /p-\alpha )}\Vert u\Vert_{B^\alpha_{p,\infty}(\Omega ,\mu )} \\
&\lesssim R^\sigma \Vert u\Vert_{B^\alpha_{p,\infty}(\Omega ,\mu )}
\end{align*}
and a similar estimate holds replacing $-N_R$ by $1$ (getting then $2^\sigma$ instead of $R^\sigma$). Replacing the use of the Hölder inequality in the proof of the case $q<\infty$ with these estimates yields
$$\lim_{R\to 0^+}\|T_Ru-T_ru\|_{L^p(\partial\Omega,\nu)} =0\quad\mbox{ and }\quad F\in L^1(\Omega ,\mu ),$$
while the rest of the proof is exactly as when $q<\infty$.
\end{proof}

\subsection{Energy bounds.}
Equipped with the chain of balls given by Lemma~\ref{lem:uniformcover}, we are now able to prove the energy bounds for the trace operator defined in Proposition~\ref{prop:tracelimit}.

\begin{prop}\label{thm:TraceEnergy}
Let $(\Omega,d,\mu)$ be a locally compact, non-complete metric measure space, with $\mu$ a doubling measure, such that $\Omega$ is an $A$-uniform domain in its completion $\overline\Omega$.  Let $\partial\Omega:=\overline\Omega\setminus\Omega$, the boundary of $\Omega$, be equipped with a Borel measure $\nu$ which is codimension $\theta$ Ahlfors regular with respect to $\mu$.  Let $1\le p<\infty$, $1\leq q\leq\infty $ and let $0<\alpha<1$ be such that $0<\alpha-\theta/p<1$.  Then there exists a bounded linear trace operator 
\[
T:HB^\alpha_{p,q}(\Omega,\mu)\to HB^{\alpha-\theta/p}_{p,q}(\partial\Omega,\nu).
\]
That is, there exists a constant $C\ge 1$, depending only on $\alpha$, $p$, $q$,  $\theta$, $C_\mu$, $C_\theta$, and $A$, such that 
\[
\| Tu\|_{HB^{\alpha-\theta/p}_{p,q}(\partial\Omega,\nu)}\le C\|u\|_{HB^{\alpha}_{p,q}(\Omega,\mu)}
\]
for all $u\in HB^{\alpha}_{p,q}(\Omega,\mu)$. 
\end{prop}

\begin{proof}
As before, we first consider the case $1\le q<\infty$.  Consider the operator $T:HB^\alpha_{p,q}(\Omega ,\mu)\rightarrow L^p(\partial\Omega ,\nu )$ given by Proposition \ref{prop:tracelimit}. Then for $\nu$-a.e. $z\in \partial\Omega$ we have
\begin{equation} \label{eq:tracelimit}
Tu(z):=\lim_{r\to0^+}\fint_{B(z,r) }u\,d\mu.
\end{equation}
Let $z,w\in\partial\Omega $ such that \eqref{eq:tracelimit} holds for both of them, and consider the corresponding chain of balls $
\{ B_k:=B(x_k,r_k) \}_{k\in\Z}$ provided by Lemma \ref{lem:uniformcover}. For the reader's convenience we recall here the properties of the collection $\{ B_k\}_k$:
\begin{enumerate}
    \item[(i)]For each $k\in\Z$, we have that $8B_k\subset\Omega$,
    \item[(ii)]$\lim_{k\to\infty}x_k=z$, $\lim_{k\to-\infty}x_k=w$,
    \item[(iii)] for each $k\in\Z$, we have that $2^{-1}B_k\cap 2^{-1}B_{k+1}\ne\varnothing$,
    \item[(iv)] there exists $N\ge 1$, such that for each $k\ge 0$, we have that $2^{-1}r_k \leq r_{k+1}\leq r_k$, and for each $x\in B_k$, we have that \[r_k\simeq d(x,\partial\Omega)\simeq d(x,z)\lesssim 2^{-|k|/N}d(z,w),\]
    \item[(v)] there exists $N\ge 1$ such that for each $k<0$, we have that $2^{-1}r_k \leq r_{k-1}\leq r_k$, and for each $x\in B_k$, we have that \[r_k\simeq d(x,\partial\Omega)\simeq d(x,w)\lesssim 2^{-|k|/N}d(z,w),\]
    \item[(vi)] there exists $C\ge 1$ such that 
    \[
    \sum_{k}\chi_{4B_k}\le C,
    \]
    \end{enumerate}
    Here the constants $C$ and $N$ and the comparison constants depend only on $A$ and $C_\mu$.
Note that when $k\ge 0$, it follows from these properties that $2B_k\supset B_{k+1}$, and  $2B_k\supset B_{k-1}$ when $k<0$.  Since $$\lim_{k\to\infty}|Tu(z)-u_{B_k}|=0=\lim_{k\to-\infty}|Tu(w)-u_{B_k}|,$$
it then follows from induction and doubling that 
\begin{align*}
|Tu&(z)-Tu(w)|^p\le\left(\sum_{k\ge 0}|u_{B_{k+1}}-u_{B_k}|+\sum_{k<0}|u_{B_k}-u_{B_{k-1}}|\right)^p\\
	&\lesssim\left(\sum_{k\ge 0}\fint_{2B_k}\fint_{2B_k}|u(x)-u(y)|d\mu(y)d\mu(x)\right)^p+ \left(\sum_{k<0}\fint_{2B_k}\fint_{2B_k}|u(x)-u(y)|d\mu(y)d\mu(x)\right)^p.
\end{align*}
Since $\alpha p-\theta>0$ by hypothesis, we choose $\beta>0$ sufficiently small so that $\alpha p-\theta-\beta p>0$.  It then follows from H\"older's inequality, property $(iv)$ of the chain of balls, and doubling that  
\begin{align*}
\Bigg(\sum_{k\ge 0}&\fint_{2B_k}\fint_{2B_k}|u(x)-u(y)|d\mu(y)d\mu(x)\Bigg)^p
	=\left(\sum_{k\ge 0}\frac{r_k^\beta}{r_k^\beta}\fint_{2B_k}\fint_{2B_k}|u(x)-u(y)|d\mu(y)d\mu(x)\right)^p \\
    &\lesssim \left(\sum_{k\ge 0}\frac{(2^{-k/N}d(z,w))^\beta}{d(x,\partial\Omega )^\beta}\fint_{2B_k}\fint_{2B_k}|u(x)-u(y)|d\mu(y)d\mu(x)\right)^p\\
	&\le\left(\sum_{k\ge 0}(2^{-k/N}d(z,w))^{\beta p'}\right)^{p/p'}\sum_{k\ge 0}d(x,\partial\Omega )^{-\beta p}\fint_{2B_k}\fint_{2B_k}|u(x)-u(y)|^p d\mu(y)d\mu(x)\\
	&\lesssim d(z,w)^{\beta p}\sum_{k\ge 0}\int_{2B_k}\int_{2B_k}\frac{|u(x)-u(y)|^p}{d(x,\partial\Omega)^{\beta p}\mu(B(x, d(x,\partial\Omega)))^2}d\mu(y)d\mu(x)\\
	&\lesssim d(z,w)^{\beta p}\int_{C_{z,w}^1}\int_{B(x,d(x,\partial\Omega))}\frac{|u(x)-u(y)|^p}{d(x,\partial\Omega)^{\beta p}\mu(B(x, d(x,\partial\Omega)))^2}d\mu(y)d\mu(x).
\end{align*}
Here, we use the notation $C_{z,w}^1:=\bigcup_{k\ge 0}2B_k$.  In the above, we have also used property $(iv)$, $(i)$, and $(vi)$ of the chain of balls to obtain the last two inequalities.   Similarly, letting $C_{z,w}^2:=\bigcup_{k<0} 2B_k$, we have that 
\begin{align*}
\Bigg(\sum_{k<0}&\fint_{2B_k}\fint_{2B_k}|u(x)-u(y)|d\mu(y)d\mu(x)\Bigg)^p\\
	&\lesssim d(z,w)^{\beta p}\int_{C_{z,w}^2}\int_{B(x,d(x,\partial\Omega))}\frac{|u(x)-u(y)|^p}{d(x,\partial\Omega)^{\beta p}\mu(B(x, d(x,\partial\Omega))^2}d\mu(y)d\mu(x).
\end{align*}
Hence, using \eqref{eq:SumForm}, it follows that 
\begin{align}\label{eq:trace I+II}
\|Tu\|^q_{HB^{\alpha-\theta/p}_{p,q}(\partial\Omega,\nu)}&\simeq\sum_{i\in\Z}2^{-i(\alpha-\theta/p)q}\left(\int_{\partial \Omega}\fint_{B(z,2^i)}|Tu(z)-Tu(w)|^pd\nu(w)d\nu(w)\right)^{q/p}\nonumber\\
	&\lesssim I+II,
\end{align}
where
\begin{align*}
I:=\sum_{i\in\Z}2^{-i(\alpha-\theta/p)q}\left(\int_{\partial\Omega}\int_{B(z,2^i)}\int_{C_{z,w}^1}\int_{B(x,d(x,\partial\Omega))}\frac{|u(x)-u(y)|^pd(z,w)^{\beta p}d\mu(y)d\mu(x)d\nu(w)d\nu(z)}{d(x,\partial\Omega)^{\beta p}\mu(B(x,d(x,\partial\Omega)))^2\nu(B(z,2^i))}\right)^{q/p}
\end{align*}
and 
\begin{align*}
II:=\sum_{i\in\Z}2^{-i(\alpha-\theta/p)q}\left(\int_{\partial\Omega}\int_{B(z,2^i)}\int_{C_{z,w}^2}\int_{B(x,d(x,\partial\Omega))}\frac{|u(x)-u(y)|^pd(z,w)^{\beta p}d\mu(y)d\mu(x)d\nu(w)d\nu(z)}{d(x,\partial\Omega)^{\beta p}\mu(B(x,d(x,\partial\Omega)))^2\nu(B(z,2^i))}\right)^{q/p}.
\end{align*}
Here the comparison constants depend only on $\alpha$, $p$, $q$, $C_\mu$, and $A$.

We first estimate $I$.  For each $j\in\Z$, let 
\[
\Omega_j:=\{x\in\Omega:2^{j-1}\le d(x,\partial\Omega)<2^{j}\}.
\]
For each $i\in\Z$, we note that if $z\in\partial\Omega$ and $w\in B(z,2^i)\cap\partial\Omega$, then by property $(iv)$ of the chain of balls, we have that $C_{z,w}^1\subset\bigcup_{j=-\infty}^{i+C}\Omega_j$ for some constant $C$ depending only on $A$ and $C_\mu$.  Hence, from this fact and Tonelli's theorem, we have that 
\begin{align*}
    \int_{\partial\Omega}&\int_{B(z,2^i)}\int_{C_{z,w}^1}\int_{B(x,d(x,\partial\Omega)}\frac{|u(x)-u(y)|^pd(z,w)^{\beta p}d\mu(y)d\mu(x)d\nu(w)d\nu(z)}{d(x,\partial\Omega)^{\beta p}\mu(B(x,d(x,\partial\Omega)))^2\nu(B(z,2^i))}\\
    &\lesssim\sum_{j=-\infty}^{i+C}\int_{\partial\Omega}\int_{B(z,2^i)}\int_{C_{z,w}^1\cap\Omega_j}\int_{B(x,2^j)}\frac{|u(x)-u(y)|^pd(z,w)^{\beta p}d\mu(y)d\mu(x)d\nu(w)d\nu(z)}{2^{j\beta p}\mu(B(x,2^j))^2\nu(B(z,2^i))}\\
    &=\sum_{j=-\infty}^{i+C}\int_{\partial\Omega}\int_{B(z,2^i)}\int_{\Omega}\int_{B(x,2^j)}\frac{|u(x)-u(y)|^pd(z,w)^{\beta p}\chi_{C_{z,w}^1\cap\Omega_j}(x)d\mu(y)d\mu(x)}{2^{j\beta p}\mu(B(x,2^j))^2\nu(B(z,2^i))}d\nu(w)d\nu(z)\\
    &=\sum_{j=-\infty}^{i+C}2^{-j\beta p}\int_{\Omega}\fint_{B(x,2^j)}\frac{|u(x)-u(y)|^p}{\mu(B(x,2^j))}\int_{\partial\Omega}\int_{B(z,2^i)}\frac{d(z,w)^{\beta p}\chi_{C_{z,w}^1\cap\Omega_j}(x)}{\nu(B(z,2^i))}d\nu(w)d\nu(z)d\mu(y)d\mu(x).
\end{align*}
We note that if $x\in C_{z,w}^1\cap\Omega_j$, then by property $(iv)$ of the chain of balls, there exists $C'\ge 1$, depending only on $A$ and $C_\mu$, such that $z\in\partial\Omega\cap B(x,C'2^j)$. Thus, using the $\theta$-codimensional relationship between $\nu$ and $\mu$, we have that 
\begin{align*}
    \int_{\partial\Omega}\int_{B(z,2^i)}\frac{d(z,w)^{\beta p}\chi_{C_{z,w}^1\cap\Omega_j}(x)}{\nu(B(z,2^i))}d\nu(w)d\nu(z)&=\int_{\partial\Omega\cap B(x,C'2^j)}\int_{B(z,2^i)}\frac{d(z,w)^{\beta p}}{\nu(B(z,2^i))}d\nu(w)d\nu(z)\\
    &\le 2^{i\beta p}\nu(B(x,C'2^j)\cap\partial\Omega)\\
    &\lesssim 2^{i\beta p-j\theta}\mu(B(x,2^j)).
\end{align*}
Substituting this into the previous expression, we then obtain the following estimate for $I$:
\begin{align*}
    I&\lesssim\sum_{i\in\Z}\left(\sum_{j=-\infty}^{i+C}2^{-i(\alpha-\theta/p)p-j\beta p+i\beta p-j\theta}\int_\Omega\fint_{B(x,2^j)}|u(x)-u(y)|^pd\mu(y)d\mu(x)\right)^{q/p}.
\end{align*}

For $i\in\Z$ and $j\in\Z$ with $j\le i+C$, we have that 
\begin{align*}
    2^{-i(\alpha-\theta/p)p-j\beta p+i\beta p-j\theta}&=\left(2^{\alpha p-\theta-\beta p}\right)^C\left(2^{\alpha p-\theta-\beta p}\right)^{-|j-(i+C)|}2^{-j\alpha p}
\end{align*}
Setting $a:=2^{\alpha p-\theta-\beta p}$, we have from our choice of $\beta$ that $1<a<\infty$.  Setting $b:=q/p$ and 
\[
c_j:=2^{-j\alpha p}\int_\Omega\fint_{B(x,2^j)}|u(x)-u(y)|^pd\mu(y)d\mu(x),
\]
we then have that 
\[
I\lesssim a^{Cb}\sum_{i\in\Z}\left(\sum_{j=-\infty}^{i+C}a^{-|j-(i+C)|}c_j\right)^b\le a^{Cb}\sum_{i\in\Z}\left(\sum_{j\in\Z}a^{-|j-i|}c_j\right)^b.
\]
It then follows from Lemma~\ref{lem:Sums} and \eqref{eq:SumForm} that 
\begin{align}\label{eq:I estimate}
    I\lesssim\sum_{j\in\Z}c_j^b=\sum_{j\in\Z}2^{-j\alpha q}\left(\int_\Omega\fint_{B(x,2^j)}|u(x)-u(y)|^pd\mu(y)d\mu(x)\right)^{q/p}\simeq \|u\|_{HB^\alpha_{p,q}(\Omega,\mu)}^q,
\end{align}
with comparison constants depending only on $\alpha$, $p$, $q$, $\theta$, $C_\mu$, $C_\theta$ and $A$. 

By a similar argument, with the roles of $w$ and $z$ reversed, we also obtain the estimate  
\[
II\lesssim\|u\|^q_{HB^\alpha_{p,q}(\Omega,\mu)},
\]
with the same dependencies for the comparison constant.  Combining these two estimates with \eqref{eq:trace I+II}, we have that 
\[
\|Tu\|_{B^{\alpha-\theta/p}_{p,q}(\partial\Omega,\nu)}\lesssim\|u\|_{B^\alpha_{p,q}(\Omega,\mu)}.
\]

For the case $q=\infty$ set, for each $i\in\Z$,
\begin{align*}
    I_i&:=\left(\sum_{j=-\infty}^{i+C}2^{-i(\alpha-\theta/p)p-j\beta p+i\beta p-j\theta}\int_\Omega\fint_{B(x,2^j)}|u(x)-u(y)|^pd\mu(y)d\mu(x)\right)^{1/p}.
\end{align*}
By the definition of the $HB^\alpha_{p,\infty}$-seminorm and the above arguments, it suffices to prove $I_i\lesssim \Vert u\Vert_{B^\alpha_{p,\infty}(\Omega ,\mu )}$ for each $i\in\Z$.  Defining $a$ and $c_j$ as before, it is clear that for any $j\in \Z$ one has $c_j\leq \Vert u\Vert^p_{B^\alpha_{p,\infty}(\Omega ,\mu )}$ and so estimating $I_i$ analogously as it was done for $I$, we have
$$I_i\lesssim a^{C/p}\left(\sum_{j=-\infty}^{i+C}a^{-|j-(i+C)|}c_j\right)^{\frac{1}{p}}\leq  a^{C/p} \Vert u\Vert_{B^\alpha_{p,\infty}(\Omega ,\mu )}\left(\sum_{j=-\infty}^{i+C}a^{-|j-(i+C)|}\right)^{\frac{1}{p}}\lesssim  \Vert u\Vert_{B^\alpha_{p,\infty}(\Omega ,\mu )}.\qedhere$$
\end{proof}

Theorem~\ref{thm:Trace} and Theorem~\ref{thm:Intro Inhom Trace} are now proved by combining Proposition~\ref{prop:tracelimit} and Proposition~\ref{thm:TraceEnergy}.

\section{Proof of Extension Theorems~\ref{thm:Extension} and \ref{thm:InhomExtension}}\label{sec:Extension}
In this section, we assume that $(\Omega,d,\mu)$ is a locally compact, non-complete metric measure space, with $\mu$ a doubling measure. We extend $\mu$ to $\overline\Omega$ by the zero extension, and note that the metric measure space $(\overline\Omega,d,\mu)$ is still doubling with the same constant. We equip  $\partial\Omega:=\overline\Omega\setminus\Omega$, the boundary of $\Omega$, with a Borel measure $\nu$ which is codimension $\theta$ Ahlfors regular with respect to $\mu$ for some $\theta>0$. Recall that this implies that $\nu$ is doubling, see Remark~\ref{rem:doublingboundary}.  Unlike the previous section, we do not assume that $\Omega$ is a uniform domain.

We first prove Theorem~\ref{thm:Extension}, constructing an extension operator for the homogeneous Besov spaces using standard Whitney extension techniques.
Let $\{ B_{i,j}\}_{i,j}$ be the Whitney cover of $\Omega$ given by Lemma \ref{thm:Whitney}, and for each $B_{i,j}=B(p_{i,j},r_{i,j})$, we consider its corresponding ``shadow'' on the boundary, given by 
\begin{equation}\label{eq:Uij}
U_{i,j}:=B(q_{i,j},r_{i,j})\cap\partial\Omega ,
\end{equation}
where $q_{i,j}\in\partial\Omega$ is a closest point to $p_{i,j}$.  By Lemma~\ref{thm:Whitney} (v), the collection $\{B_{i,j}\}_j$ has bounded overlap for each $i$, and so this implies bounded overlap of the collection $\{U_{i,j}\}_j$, as shown by the following lemma:

\begin{lem}\label{lem:bounedoverlapUij}
Fix $i\in\Z$.  For each $K\ge 1$, there exists $C_K'\ge 1$ such that 
\[
\sum_j\chi_{KU_{i,j}}\le C_K',
\]
where $KU_{i,j}=B(q_{i,j},Kr_{i,j})\cap\partial\Omega$.
\end{lem}

\begin{proof}
    If $z\in KU_{i,j}$, then we have that 
    \[
    d(z,p_{i,j})\le d(z,q_{i,j})+d(q_{i,j},p_{i,j})<Kr_{i,j}+8r_{i,j},
    \]
    and so $KU_{i,j}\subset (K+8)B_{i,j}$.  Therefore by Lemma~\ref{thm:Whitney} (v), it follows that 
    \[
    \sum_j\chi_{KU_{i,j}}\le\sum_{j}\chi_{(K+8)B_{i,j}}\le C_{K+8}=:C_K'.\qedhere
    \]
\end{proof}

We now prove Theorem~\ref{thm:Extension}:

\begin{proof}[Proof of Theorem~\ref{thm:Extension}]
We first consider the case $1\le q<\infty$.  
Let $f\in HB^{\alpha}_{p,q}(\partial\Omega,\nu)$, and let $\mathcal{W}_\Omega=\{B(p_{i,j},r_{i,j})=:B_{i,j}\}$ and $\{\phii_{i,j}\}$ be the Whitney cover and partition of unity given by Lemma~\ref{thm:Whitney}.  For each $(i,j)$  let $U_{i,j}$ as in \eqref{eq:Uij} 
and $a_{i,j}:=\fint_{U_{i,j}}fd\nu.$  Then, for each $x\in\Omega,$ let 
\begin{equation}\label{eq:homExt}
\wtil Ef(x):=\sum_{i,j}a_{i,j}\phii_{i,j}(x).
\end{equation}
By this construction, the map $f\mapsto \wtil Ef$ is linear.  By \eqref{eq:SumForm}, it follows that 
\begin{align}\label{eq:First}
\|\wtil Ef\| ^q_{HB^{\alpha+\theta /p}_{p,q}(\Omega,\mu)}&\simeq\sum_{k\in\Z}2^{-k(\alpha+\theta/p)q}\left(\int_{\Omega}\fint_{B(x,2^k)}|\wtil Ef(x)-\wtil Ef(y)|^pd\mu(y)d\mu(x)\right)^{q/p}\nonumber\\
	&\le\sum_{k\in\Z}2^{-k(\alpha+\theta/p)q}\left(\sum_{i,j}\int_{B_{i,j}}\fint_{B(x,2^k)}|\wtil Ef(x)-\wtil Ef(y)|^pd\mu(y)d\mu(x)\right)^{q/p}.
\end{align}
Fixing $k\in\Z$, we then have
\begin{align}\label{eq:Inside parenth I+II}
    \sum_{i,j}\int_{B_{i,j}}\fint_{B(x,2^k) }&|\wtil Ef(x)-\wtil Ef(y)|^pd\mu(y)d\mu(x)\nonumber\\
        &=\sum_{\substack{i,j\st\\ i\ge k+1}}\int_{B_{i,j}}\fint_{B(x,2^k) }|\wtil Ef(x)-\wtil Ef(y)|^pd\mu(y)d\mu(x)\nonumber\\
        &\qquad+\sum_{\substack{i,j\st\\i<k+1}}\int_{B_{i,j}}\fint_{B(x,2^k) }|\wtil Ef(x)-\wtil Ef(y)|^pd\mu(y)d\mu(x)=:I+II.
\end{align}

To estimate $I$, let $(i,j)$ be such that $i\ge k+1.$  Then for $x\in B_{i,j}$ and $y\in B(x,2^k)$, we have that $B(x,2^k)\subset 2B_{i,j}$, and so by the partition of unity, it follows that 
\begin{align*}
|\wtil Ef(x)-\wtil Ef(y)|&=\left|\sum_{l,m} a_{l,m}\varphi_{l,m}(x)-\sum_{l,m}a_{l,m}\varphi_{l,m}   (y)\right|\\
    &= \left| \sum_{l,m}a_{l,m}\varphi_{l,m}(x)-a_{i,j}\sum_{l,m}\varphi_{l,m}(x)+a_{i,j}\sum_{l,m}\varphi_{l,m} (y)-\sum_{l,m}a_{l,m}\varphi_{l,m}(y)\right| \\
    &\le \sum_{l,m}|a_{l,m}-a_{i,j}||\phii_{l,m}(x)-\phii_{l,m}(y)|=\sum_{\substack{l,m\st\\2B_{i,j}\cap 2B_{l,m}\ne\varnothing}}|a_{l,m}-a_{i,j}||\phii_{l,m}(x)-\phii_{l,m}(y)|.
\end{align*}
The last equality follows since $\phii_{l,m}$ is supported in $2B_{l,m}$.  

If $(l,m)$ is such that $2B_{i,j}\cap 2B_{l,m}\ne\varnothing$, then since $r_{i,j}=\frac{1}{8}d(p_{i,j},\partial\Omega)$, it follows that $2^{-1}r_{i,j}\le r_{l,m}\le 2r_{i,j}$, and so $i-2\le l\le i+2$ due to (iii) in Lemma \ref{thm:Whitney}.  By the triangle inequality, we then have that 
\begin{equation}\label{eq:U^*}
U_{i,j}^*:=B(q_{i,j},32r_{i,j})\cap\partial\Omega\supset U_{l,m}.
\end{equation}
Furthermore, it follows that $B_{i,j}\subset 8B_{l,m}$, and so by Lemma~\ref{thm:Whitney} (v), it follows that there are at most $C\ge 1$ indices $(l,m)$ such that $2B_{l,m}\cap 2B_{i,j}\ne\varnothing$, with $C$ depending only on $C_\mu$.  Finally, since $r_{i,j}\simeq r_{l,m}$, we have that $\phii_{l,m}$ is $C/r_{i,j}$-Lipschitz, with $C$ depending only on the doubling constant. Using these facts, the doubling property of $\nu$, and Jensen's inequality, it follows that 
\begin{align*}
|\wtil Ef(x)-\wtil Ef(y)|^p&\lesssim \left(\frac{d(x,y)}{r_{i,j}}\sum_{\substack{l,m\st\\2B_{i,j}\cap 2B_{l,m}\ne\varnothing}}|a_{l,m}-a_{i,j}|\right)^p\\
	&\lesssim \left( \frac{d(x,y)}{r_{i,j}}\sum_{\substack{l,m\st\\2B_{i,j}\cap 2B_{l,m}\ne\varnothing}}\fint_{U_{l,m}}\fint_{U_{i,j}}|f(w)-f(z)|d\nu(w)d\nu(z)\right)^p \\
	&\lesssim \left( \frac{d(x,y)}{r_{i,j}}\fint_{U_{i,j}^*}\fint_{U_{i,j}^*}|f(w)-f(z)|d\nu(w)d\nu(z)\right)^p\\
    &\le\frac{d(x,y)^p}{r^p_{i,j}}\fint_{U_{i,j}^*}\fint_{U_{i,j}^*}|f(w)-f(z)|^pd\nu(w)d\nu(z),
\end{align*}
with comparison constant depending only on  $\theta$, $C_\mu$ and $C_\theta$.  Hence, we have that 
\begin{align*}
\int_{B_{i,j}}\fint_{B(x,2^k)}|\wtil Ef(x)&-\wtil Ef(y)|^pd\mu(y)d\mu(x)\\
	&\lesssim\frac{1}{r^p_{i,j}}\fint_{U_{i,j}^*}\fint_{U_{i,j}^*}|f(w)-f(z)|^pd\nu(w)d\nu(z)\int_{B_{i,j}}\fint_{B(x,2^k)}d(x,y)^pd\mu(y)d\mu(x)\\	
 &\le\frac{1}{r^p_{i,j}}\fint_{U_{i,j}^*}\fint_{U_{i,j}^*}|f(w)-f(z)|^pd\nu(w)d\nu(z)\int_{B_{i,j}}\fint_{B(x,2^k)}2^{pk}d\mu(y)d\mu(x)\\	
	&\lesssim 2^{p(k-i)}\mu(B_{i,j})\fint_{U_{i,j}^*}\fint_{U_{i,j}^*}|f(w)-f(z)|^pd\nu(w)d\nu(z)\\
	&\le2^{p(k-i)}\frac{\mu(B(q_{i,j},32r_{i,j}))}{\nu(U_{i,j}^*)}\int_{U_{i,j}^*}\fint_{U_{i,j}^*}|f(w)-f(z)|^pd\nu(w)d\nu(z)\\
	&\lesssim 2^{p(k-i)} r_{i,j}^\theta\int_{U_{i,j}^*}\fint_{U_{i,j}^*}|f(w)-f(z)|^pd\nu(w)d\nu(z)\\
    &\le2^{p(k-i)+i\theta}\int_{U_{i,j}^*}\fint_{U_{i,j}^*}|f(w)-f(z)|^pd\nu(w)d\nu(z),
\end{align*}
with the comparison constant depending only on $\theta$, $C_\mu$, and $C_\theta$. By Lemma~\ref{lem:bounedoverlapUij}, we then have that 
\begin{align}\label{eq:Inside Parenth I (1)}
    I\lesssim\sum_{i\ge k+1}2^{p(k-i)+i\theta}\int_{\partial\Omega}\fint_{B(z,2^{i+6})\cap\partial\Omega}|f(w)-f(z)|^pd\nu(w)d\nu(z),
\end{align}

To estimate $II$ in \eqref{eq:Inside parenth I+II}, we consider $i,j$ such that $i<k+1$.

For each $x\in B_{i,j}$, we have that $B(x,2^k)\subset B(p_{i,j},2^{k+1})$.  Letting
\[
I_{i,j,k}:=\{l,m: B_{l,m}\cap B(p_{i,j},2^{k+1})\ne\varnothing\},
\]
we have that 
\begin{align*}
\int_{B_{i,j}}\fint_{B(x,2^k) }|\wtil Ef(x)-\wtil Ef(y)|^pd\mu(y)d\mu(x)\nonumber&\le\sum_{l,m\in I_{i,j,k}}\int_{B_{i,j}}\int_{B_{l,m}}\frac{|\wtil Ef(x)-\wtil Ef(y)|^p}{\mu(B(x,2^k) )}d\mu(y)d\mu(x).
\end{align*}
For $x\in B_{i,j}$ and $y\in B_{l,m}$ for $(l,m)\in I_{i,j,k},$ it follows from the partition of unity, Lemma~\ref{thm:Whitney} (ii) and Jensen's inequality that 
\begin{align*}
|\wtil Ef(x)&-\wtil Ef(y)|^p \le\left( \sum_{s,t}|a_{s,t}-a_{i,j}||\phii_{s,t}(x)-\phii_{s,t}(y)|\right)^p\\
	&\le \left( \sum_{\substack{s,t\st\\x\in 2B_{s,t}}}|a_{s,t}-a_{i,j}|+\sum_{\substack{s,t\st\\y\in 2B_{s,t}}}|a_{s,t}-a_{i,j}|\right)^p\\
	&\lesssim \left( \fint_{U_{i,j}^*}\fint_{U_{i,j}^*}|f(w)-f(z)|d\nu(w)d\nu(z)+\fint_{U_{i,j}^*}\fint_{U_{l,m}^*}|f(w)-f(z)|d\nu(w)d\nu(z)\right)^p\\
 &\lesssim\fint_{U_{i,j}^*}\fint_{U_{i,j}^*}|f(w)-f(z)|^pd\nu(w)d\nu(z)+\fint_{U_{i,j}^*}\fint_{U_{l,m}^*}|f(w)-f(z)|^pd\nu(w)d\nu(z),
\end{align*}
where the comparison constant depends only on $\theta$, $p$, $C_\mu$ and $C_\theta$.  Substituting this into the above expression, we have that 
\begin{align}\label{eq:Second}
\int_{B_{i,j}}&\fint_{B(x,2^k) }|\wtil Ef(x)-\wtil Ef(y)|^pd\mu(y)d\mu(x)\nonumber\\
	&\lesssim \fint_{U_{i,j}^*}\fint_{U_{i,j}^*}|f(w)-f(z)|^pd\nu(w)d\nu(z)\sum_{l,m\in I_{i,j,k}}\int_{B_{i,j}}\int_{B_{l,m}}\frac{1}{\mu(B(x,2^k) )}d\mu(y)d\mu(x)\nonumber\\
	&+\sum_{l,m\in I_{i,j,k}}\left(\fint_{U_{i,j}^*}\fint_{U_{l,m}^*}|f(w)-f(z)|^pd\nu(w)d\nu(z)\int_{B_{i,j}}\int_{B_{l,m}}\frac{1}{\mu(B(x,2^k) )}d\mu(y)d\mu(x)\right)\nonumber\\
	&=:III+IV.
\end{align}

For $(l,m)\in I_{i,j,k}$, we have by Lemma~\ref{thm:Whitney} and the assumption that $i<k+1$ that
\begin{align*}
8r_{l,m}=d(p_{l,m},\partial\Omega)\le r_{l,m}+2^{k+1}+d(p_{i,j},\partial\Omega)\le r_{l,m}+2^{k+1}+8r_{i,j}\le r_{l,m}+10\cdot 2^{k},
\end{align*}
and so it follows that 
\[
3\cdot 2^{l}\le 10\cdot 2^k.
\]
Hence, for $(l,m)\in I_{i,j,k}$, we have that $l\le k+2$, and so if $x\in B_{i,j}$ and $y\in B_{l,m}$ for $(l,m)\in I_{i,j,k}$, it follows that 
\[
d(x,y)\le 2r_{l,m}+r_{i,j}+2^{k+1}\le 2^{k+4}.
\]
Therefore, by bounded overlap, we can sum over $(l,m)\in I_{i,j,k}$ to estimate $III$ as follows: 

\[
III\lesssim\fint_{U_{i,j}^*}\fint_{U_{i,j}^*}|f(w)-f(z)|^pd\nu(w)d\nu(z)\int_{B_{i,j}}\int_{B(x,2^{k+4}) }\frac{1}{\mu(B(x,2^k) )}d\mu(y)d\mu(x).
\]
By doubling of $\mu$ and the codimensionality of $\nu$, it follows that 
\begin{align}\label{eq:III}
III&\lesssim\frac{\mu(B_{i,j})}{\nu(U_{i,j}^*)}\int_{U_{i,j}^*}\fint_{U_{i,j}^*}|f(w)-f(z)|^pd\nu(w)d\nu(z)\nonumber\\
    &\simeq\frac{\mu(B(q_{i,j}, 32r_{i,j}))}{\nu(U_{i,j}^*)}\int_{U_{i,j}^*}\fint_{U_{i,j}^*}|f(w)-f(z)|^pd\nu(w)d\nu(z)\nonumber\\
    &\simeq  r_{i,j}^\theta\int_{U_{i,j}^*}\fint_{U_{i,j}^*}|f(w)-f(z)|^pd\nu(w)d\nu(z)\le 2^{i\theta}\int_{U_{i,j}^*}\fint_{U_{i,j}^*}|f(w)-f(z)|^pd\nu(w)d\nu(z)
\end{align}
with comparison constant depending only on $\theta$, $C_\mu$ and $C_\theta$.

Now, for $x\in B_{i,j},$ we have that 
$\mu(B(x,2^k) )\simeq\mu(B(q_{i,j},2^k) )$ by the doubling condition of $\mu$.  Using this, we have that 
\begin{align*}
IV&\simeq \frac{1}{\mu(B(q_{i,j},2^k) )}\sum_{l,m\in I_{i,j,k}}\frac{\mu(B_{i,j})\mu(B_{l,m})}{\nu(U_{i,j}^*)\nu(U_{l,m}^*)}\int_{U_{i,j}^*}\int_{U_{l,m}^*}|f(w)-f(z)|^pd\nu(w)d\nu(z)\nonumber\\
	&\lesssim\frac{2^{i\theta}}{\mu(B(q_{i,j},2^k) )}\sum_{l,m\in I_{i,j,k}}2^{l\theta}\int_{U_{i,j}^*}\int_{U_{l,m}^*}|f(w)-f(z)|^pd\nu(w)d\nu(z),\nonumber\\
\end{align*}
with comparison constant again depending only on  $\theta$, $C_\mu$ and $C_\theta$.

Recall that if $(l,m)\in I_{i,j,k}$, by definition we have $d(p_{i,j},p_{l,m})\leq 2^{k+1}+r_{l,m}$. Moreover, we also saw that $l\le k+2$, and so we have that

\[
d(q_{i,j},q_{l,m})\leq 8(r_{i,j}+r_{l,m})+d(p_{i,j},p_{l,m})<8\cdot 2^{k}+9\cdot 2^{k+2}+2^{k+1}\leq 2^{k+6},
\]
Therefore, for $z\in U^*_{i,j}$ and $w\in U^*_{l,m}$ it follows that 
$
d(z,w)\leq 32(r_{i,j}+r_{l,m})+d(q_{i,j},q_{l,m})\leq 2^{k+8}$
and so $U^*_{l,m}\subset B(z,2^{k+8})\cap\partial\Omega$.  Using this fact, Lemma~\ref{lem:bounedoverlapUij}, and the codimensionality and doubling property  of $\nu$, we have that 
\begin{align}\label{eq:IV}
IV&\lesssim \frac{2^{i\theta}}{\mu(B(q_{i,j},2^k)  )}\sum_{l=-\infty}^{k+2}\sum_{\substack{m\st\\l,m\in I_{i,j,k}}}2^{l\theta}\int_{U_{i,j}^*}\int_{U_{l,m}^*}|f(w)-f(z)|^pd\nu(w)d\nu(z)\nonumber\\
	&\lesssim \frac{2^{i\theta}}{\mu(B(q_{i,j},2^k)  )}\int_{U_{i,j}^*}\int_{B(z,2^{k+ 8})\cap\partial\Omega}|f(w)-f(z)|^pd\nu(w)d\nu(z)\sum_{l=-\infty}^{k+2}2^{l\theta}\nonumber\\
	&\simeq \frac{2^{i\theta+k\theta}}{\mu(B(q_{i,j},2^k)  )}\int_{U_{i,j}^*}\int_{B(z,2^{k+ 8})\cap\partial\Omega}|f(w)-f(z)|^pd\nu(w)d\nu(z)\nonumber\\
	&\lesssim \frac{2^{i\theta}}{\nu(B(q_{i,j},2^k)\cap\partial\Omega)}\int_{U_{i,j}^*}\int_{B(z,2^{k+ 8})\cap\partial\Omega}|f(w)-f(z)|^pd\nu(w)d\nu(z)\nonumber\\
	&\lesssim2^{i\theta}\int_{U_{i,j}^*}\fint_{B(z,2^{k+ 8})}|f(w)-f(z)|^pd\nu(w)d\nu(z),
\end{align}
with comparison constant depending only on $\theta$, $C_\mu$  and $C_\theta$.  In particular the constant depends on the convergence of $\sum_{l=-\infty}^0 2^{l\theta}$.

Combining \eqref{eq:Inside parenth I+II}, \eqref{eq:Second}, \eqref{eq:III}, and \eqref{eq:IV}, and using Lemma~\ref{lem:bounedoverlapUij}, we have that 
\begin{align}\label{eq:IIdecomp}
II\lesssim \sum_{i<k+1}&\sum_j 2^{i\theta}\int_{U_{i,j}^*}\fint_{U_{i,j}^*}|f(w)-f(z)|^pd\nu(w)d\nu(z)\nonumber\\
	&+\sum_{i<k+1} 2^{i\theta}\int_{\partial\Omega}\fint_{B(z,2^{k+ 8})\cap\partial\Omega}|f(w)-f(z)|^pd\nu(w)d\nu(z).
\end{align}
By \eqref{eq:First}, \eqref{eq:Inside parenth I+II}, \eqref{eq:Inside Parenth I (1)}, and \eqref{eq:IIdecomp}, we then have that 
\begin{align}\label{eq:V+VI+VII}
   \|\wtil Ef\|&^q_{HB^{\alpha+\theta /p}_{p,q}(\Omega,\mu)}\nonumber\\
   &\lesssim\sum_{k\in\Z}2^{-k(\alpha+\theta/p)q}\left(\sum_{i\ge k+1}2^{p(k-i)+i\theta}\int_{\partial\Omega}\fint_{B(z,2^{i+ 6})\cap\partial\Omega}|f(w)-f(z)|^pd\nu(w)d\nu(z)\right)^{q/p}\nonumber\\
   &+\sum_{k\in\Z}2^{-k(\alpha+\theta/p)q}\left(\sum_{i<k+1}\sum_j 2^{i\theta}\int_{U_{i,j}^*}\fint_{U_{i,j}^*}|f(w)-f(z)|^pd\nu(w)d\nu(z)\right)^{q/p}\nonumber\\
   &+\sum_{k\in\Z}2^{-k(\alpha+\theta/p)q}\left(\sum_{i<k+1} 2^{i\theta}\int_{\partial\Omega}\fint_{B(z,2^{k+ 8})\cap\partial\Omega}|f(w)-f(z)|^pd\nu(w)d\nu(z)\right)^{q/p}\nonumber\\
   &=:V+VI+VII.
\end{align}

We note that for $i\ge k+1$, we have that 
\begin{align*}
    2^{-k(\alpha p+\theta)+p(k-i)+i\theta}=2^{-k\alpha p+i\alpha p-i\alpha p}\left(2^{p-\theta}\right)^{-|i-k|}=\left(2^{p-\theta-\alpha p}\right)^{-|i-k|}2^{-i\alpha p}.
\end{align*}

Let $a=2^{p-\theta-\alpha p}$, $b=q/p$, and let
\[
c_i=2^{-i\alpha p}\int_{\partial\Omega}\fint_{B(z, 2^{i+ 6})\cap\partial\Omega}|f(w)-f(z)|^pd\nu(w)d\nu(z).
\]
Note that $1<a<\infty$ by our assumption that $\alpha<1-\theta/p$, we then use Lemma~\ref{lem:Sums} and \eqref{eq:SumForm} to estimate $V$ by 
\begin{align}\label{eq:V}
    V&=\sum_{k\in\Z}\left(\sum_{i\ge k+1}2^{-k(\alpha p+\theta)}2^{p(k-i)+i\theta}\int_{\partial\Omega}\fint_{B(z,2^{i+ 6})\cap\partial\Omega}|f(w)-f(z)|^pd\nu(w)d\nu(z)\right)^{q/p}\nonumber\\
    &\le\sum_{k\in\Z}\left(\sum_{i\in\Z}\left(2^{p-\theta-\alpha p}\right)^{-|i-k|}2^{-i\alpha p}\int_{\partial\Omega}\fint_{B(z,2^{i+ 6})\cap\partial\Omega}|f(w)-f(z)|^pd\nu(w)d\nu(z)\right)^{q/p}\nonumber\\
    &=\sum_{k\in\Z}\left(\sum_{i\in\Z}a^{-|i-k|}c_i\right)^b\nonumber\\
    &\lesssim \sum_{i\in\Z}c_i^b\lesssim\sum_{i\in\Z}2^{-i\alpha q}\left(\int_{\partial\Omega}\fint_{B(z,2^{i+6})\cap\partial\Omega}|f(w)-f(z)|^pd\nu(w)d\nu(z)\right)^{q/p}\simeq\|f\|_{HB^\alpha_{p,q}(\partial\Omega,\nu)}^q.
\end{align}

To estimate $VI$, we note that for $i<k+1$,
\[
2^{-k(\alpha p+\theta)+i\theta}=\left(2^{\alpha p+\theta}\right)^{-|i-k|}2^{-i\alpha p}.
\]
Setting $a=2^{\alpha p+\theta}$, $b=q/p$, and 
\[
c_i=2^{-i\alpha p}\int_{\partial\Omega}\fint_{B(z,2^{i+ 6})\cap\partial\Omega}|f(w)-f(z)|^pd\nu(w)d\nu(z),
\]
we have by Lemma~\ref{lem:bounedoverlapUij}, Lemma~\ref{lem:Sums}, and \eqref{eq:SumForm},
\begin{align}\label{eq:VI}
    VI&=\sum_{k\in\Z}\left(\sum_{i<k+1}2^{-k(\alpha p+\theta)+i\theta}\sum_j \int_{U_{i,j}^*}\fint_{U_{i,j}^*}|f(w)-f(z)|^pd\nu(w)d\nu(z)\right)^{q/p}\nonumber\\
    &\lesssim\sum_{k\in\Z}\left(\sum_{i\in\Z}\left(2^{\alpha p+\theta}\right)^{-|i-k|}2^{-i\alpha p} \int_{\partial\Omega}\fint_{B(z,2^{i+ 6})}|f(w)-f(z)|^pd\nu(w)d\nu(z)\right)^{q/p}\nonumber\\
    &=\sum_{k\in\Z}\left(\sum_{i\in\Z}a^{-|i-k|}c_i\right)^{b}\nonumber\\
    &\lesssim \sum_{i\in\Z}c_i^b\lesssim\sum_{i\in\Z}2^{-i\alpha q}\left(\int_{\partial\Omega}\fint_{B(z,2^{i+6})}|f(w)-f(z)|^pd\nu(w)d\nu(z)\right)^{q/p}\simeq\|f\|_{HB^\alpha_{p,q}(\partial\Omega,\nu)}^q.
\end{align}

Finally, since $\sum_{i<k+1}2^{i\theta}\simeq 2^{k\theta}$, it follows that 
\begin{align*}
 VII\lesssim\sum_{k\in\Z}2^{-k\alpha q}\left(\int_{\partial\Omega}\fint_{B(z,2^{k+ 8})\cap\partial\Omega}|f(w)-f(z)|^pd\nu(w)d\nu(z)\right)^{q/p}\simeq\|f\|_{HB^\alpha_{p,q}(\partial\Omega,\nu)}^q.   
\end{align*}
Notice that the estimates of $V$ and $VI$ add the dependency on $p$, $q$, $\alpha$ and $\theta$ to the comparison constant due to the use of Lemma \ref{lem:Sums} with our choice of $a$ and $b$. Combining this estimate with \eqref{eq:V}, \eqref{eq:VI}, and \eqref{eq:V+VI+VII} completes the proof for the case $1\leq q<\infty $.

For the case $q=\infty$, we follow the same arguments as above until reaching \eqref{eq:V+VI+VII}, and now instead set \begin{align*}
   &V_k:=2^{-k(\alpha+\theta/p)}\left(\sum_{i\ge k+1}2^{p(k-i)+i\theta}\int_{\partial\Omega}\fint_{B(z,2^{i+ 6})\cap\partial\Omega}|f(w)-f(z)|^pd\nu(w)d\nu(z)\right)^{1/p}\nonumber\\
   &VI_k:=2^{-k(\alpha+\theta/p)}\left(\sum_{i<k+1}\sum_j 2^{i\theta}\int_{U_{i,j}^*}\fint_{U_{i,j}^*}|f(w)-f(z)|^pd\nu(w)d\nu(z)\right)^{1/p}\nonumber\\
   &VII_k:=2^{-k(\alpha+\theta/p)}\left(\sum_{i<k+1} 2^{i\theta}\int_{\partial\Omega}\fint_{B(z,2^{k+ 8})\cap\partial\Omega}|f(w)-f(z)|^pd\nu(w)d\nu(z)\right)^{1/p}
\end{align*}
for each $k\in \Z$. We then have to prove $V_k+VI_k+VII_k\lesssim \|f\|_{HB^\alpha_{p,\infty}(\partial\Omega,\nu)} $ for every $k\in \Z$. Notice that given any $i\in\Z$ then $c_i\lesssim \|f\|_{HB^\alpha_{p,\infty}(\partial\Omega,\nu)}^p$ by definition, and thus, estimating analogously to the case $q<\infty$ we get 
$$V_k\lesssim \left(\sum_{i\in\Z}a^{-|i-k|}c_i\right)^{\frac{1}{p}} \lesssim \|f\|_{HB^\alpha_{p,\infty}(\partial\Omega,\nu)}\left(\sum_{i\in\Z}a^{-|i-k|}\right)^{\frac{1}{p}}\simeq \|f\|_{HB^\alpha_{p,\infty}(\partial\Omega,\nu)},$$
while an analogous argument yields the estimate for $VI_k$. Furthermore, 
\begin{equation*}
 VII_k\lesssim 2^{-k\alpha }\left(\int_{\partial\Omega}\fint_{B(z,2^{k+ 8})\cap\partial\Omega}|f(w)-f(z)|^pd\nu(w)d\nu(z)\right)^{1/p}\lesssim \|f\|_{HB^\alpha_{p,\infty}(\partial\Omega,\nu)}.\qedhere 
\end{equation*}
\end{proof}

We now prove Corollary~\ref{cor:inhom ext bounded domain}, showing that $\wtil E$, in fact, gives a bounded linear extension operator for the inhomogeneous spaces $B^\alpha_{p,q}(\partial\Omega\,\nu)$ when $\Omega$ is bounded:

\begin{proof}[Proof of Corollary~\ref{cor:inhom ext bounded domain}]
    For each $f\in B^\alpha_{p,q}(\partial\Omega,\nu)$, let $\wtil Ef$ be defined by \eqref{eq:homExt}. The linearity of $\wtil E$ and energy bounds \eqref{eq:hom extension energy bound} follow from Theorem~\ref{thm:Extension}. It remains to show \eqref{eq:inhom Lp bounds}.  To this end, we note that as $\Omega$ is bounded, there exists $N_\Omega\in\Z$, depending only on $\diam(\Omega)$, such that $\mathcal{W}_\Omega=\bigcup_{i=-\infty}^{N_\Omega}\{B_{i,j}\}_j$.  That is, the Whitney cover of $\Omega$ given by Lemma~\ref{thm:Whitney} contains balls of radius no greater than $2^{N_\Omega}$. Therefore, by the definition of $\wtil Ef$, we have that 
    \begin{align}\label{eq:layer estimate part 1}
        \int_\Omega|\wtil Ef|^pd\mu\le\sum_{i=\infty}^{N_\Omega}\sum_j\int_{B_{i,j}}|\wtil Ef|^pd\mu
        &\le\sum_{i=-\infty}^{N_\Omega}\sum_{j}\int_{B_{i,j}}\Bigg|\sum_{l,m}\left(\fint_{U_{l,m}}fd\nu\right)\phii_{l,m}(x)\Bigg|^pd\mu(x)\nonumber\\
    &\le\sum_{i=-\infty}^{N_\Omega}\sum_{j}\int_{B_{i,j}}\Bigg|\sum_{\substack{l,m\st\\2B_{l,m}\cap B_{i,j}\ne\varnothing}}\left(\fint_{U_{l,m}}fd\nu\right)\Bigg|^pd\mu(x).
    \end{align}
    By bounded overlap, there are at most $C$ indices $(l,m)$ such that $2B_{l,m}\cap B_{i,j}\ne\varnothing$.  Furthermore, for such indices, we have that $U_{l,m}\subset U^*_{i,j}$, where $U_{i,j}^*$ is given by \eqref{eq:U^*}, and so $\nu(U_{l,m})\simeq\nu(U^*_{i,j})$ by the doubling property of $\nu$.  From these facts, Jensen's inequality, the doubling property of $\mu$, the codimension $\theta$ Ahlfors regularity of $\nu$, and Lemma~\ref{lem:bounedoverlapUij}, we have
\begin{align}\label{eq:layer estimate part 2}
    \int_\Omega|\wtil Ef|^pd\mu&\lesssim\sum_{i=-\infty}^{N_\Omega}\sum_{j}\int_{B_{i,j}}\left(\fint_{U_{i,j}^*}|f|d\nu\right)^pd\mu(x)\nonumber\\
    &\le\sum_{i=-\infty}^{N_\Omega}\sum_j\int_{B_{i,j}}\fint_{U_{i,j}^*}|f|^pd\nu d\mu(x)\nonumber\\
	&\le\sum_{i=-\infty}^{N_\Omega}\sum_{j}\frac{\mu(B(q_{i,j},32r_{i,j}))}{\nu(U_{i,j}^*)}\int_{U_{i,j}^*}|f|^pd\nu\nonumber\\
    &\lesssim \sum_{i=-\infty}^{N_\Omega}2^{i\theta}\sum_{j}\int_{U_{i,j}^*}|f|^pd\nu\lesssim \sum_{i=-\infty}^{N_\Omega}2^{i\theta}\|f\|^p_{L^p(\partial\Omega,\nu)}\lesssim\|f\|^p_{L^p(\partial\Omega,\nu)},
\end{align}
where the comparison constant depends only on $p$, $\theta$, $C_\mu$,  $C_\theta$ and $\diam(\Omega)$.
\end{proof}

If $\Omega$ is not assumed to be bounded, then an extension result for the inhomogeneous Besov spaces can still be obtained by applying a Lipschitz cutoff function to the extension operator obtained in the previous theorem. We now prove Theorem~\ref{thm:InhomExtension}, obtaining the desired extension operator.

\begin{proof}[Proof of Theorem~\ref{thm:InhomExtension}]
    Let $f\in B^{\alpha}_{p,q}(\partial\Omega,\nu)$.  Let $\Phi:\Omega\to[0,1]$ be a $1$-Lipschitz function such that $\Phi\equiv 1$ in $\{x\in\Omega:d(x,\partial\Omega)\le 1\}$ and $\Phi\equiv 0$ in $\{x\in\Omega:d(x, \partial\Omega)>2\}$.  For $x\in\Omega$, we then set
\begin{equation*}
Ef(x)=\Phi(x)\wtil Ef(x),
\end{equation*}
where $\wtil Ef$ is given by \eqref{eq:homExt}.  

Let $\Omega_2:=\{x\in\Omega:d(x,\partial\Omega)\le 2\}$.  If $i,j$ is such that $B_{i,j}\cap\Omega_2\ne\varnothing$, then 
\[
4\cdot 2^i\le 8r_{i,j}=d(p_{i,j},\partial\Omega)\le r_{i,j}+2\le 2^i+2,
\]
and so it follows that $2^i\le 2/3$.  Hence $i<0$.  Here $\mathcal{W}_\Omega=\{B_{i,j}\}_{i,j}$ is the Whitney cover used in Theorem~\ref{thm:Extension}.  Since $0\le\Phi\le 1$, we then have that 
\begin{align}\label{eq:i<0}
    \|Ef\|^p_{L^p(\Omega,\mu)}\le\int_{\Omega_2}|\wtil Ef|^pd\mu\le\sum_{i=-\infty}^0\sum_j\int_{B_{i,j}}|\wtil Ef|^pd\mu.
\end{align}
For $i<0$, we see from estimates \eqref{eq:layer estimate part 1} and \eqref{eq:layer estimate part 2} that 
\begin{align*}
\sum_j\int_{B_{i,j}}|\wtil Ef|^pd\mu\lesssim 2^{i\theta}\|f\|^p_{L^p(\partial\Omega,\nu)}.
\end{align*}

Substituting this estimate into \eqref{eq:i<0}, we obtain
\begin{equation}\label{eq:ExtL^pBound}
\|Ef\|_{L^p(\Omega,\mu)}\le\|\wtil Ef\|_{L^p(\Omega_2,\mu)}\le C\|f\|_{L^p(\partial\Omega,\nu)},
\end{equation}
where $C\ge 1$ depends only on $p$, $\theta$, $C_\mu$  and $C_\theta$.

It remains to estimate the Besov energy of $Ef$.  Let $1\le q<\infty$ and recall by \eqref{eq:SumFormLp} and \eqref{eq:ExtL^pBound} that
\begin{align}\label{eq:InhomExt Energy bound}
\|Ef\|^q_{HB^{\alpha+\theta/p}_{p,q}(\Omega,\mu)}&\lesssim\|Ef\|_{L^p(\Omega,\mu)}^q+\sum_{i\in\N}2^{i(\alpha+\theta/p)q}\left(\int_{\Omega}\fint_{B(x,2^{-i}) }|Ef(x)-Ef(y)|^pd\mu(y)d\mu(x)\right)^{q/p}\nonumber\\
&\lesssim\|f\|_{L^p(\partial\Omega,\nu)}^q+I+II,
\end{align}
where
\[
I:=\sum_{i\in\N}2^{i(\alpha+\theta/p)q}\left(\int_{\Omega_2}\fint_{B(x,2^{-i}) }|Ef(x)-Ef(y)|^pd\mu(y)d\mu(x)\right)^{q/p}
\]
and
\[
II:=\sum_{i\in\N}2^{i(\alpha+\theta/p)q}\left(\int_{\Omega\setminus\Omega_2}\int_{B(x,2^{-i})\cap \Omega_2}\frac{|Ef(x)-Ef(y)|^p}{\mu(B(x,2^{-i}) )}d\mu(y)d\mu(x)\right)^{q/p}.
\]
By the triangle inequality, we have
\[
|Ef(x)-Ef(y)|\le|\Phi(y)||\wtil Ef(x)-\wtil Ef(y)|+|\wtil Ef(x)||\Phi(x)-\Phi(y)|,
\]
and since $\Phi$ is $1$-Lipschitz and $0\le\Phi\le 1$, it follows that 
\begin{align*}
    I&\lesssim \sum_{i\in\N}2^{i(\alpha+\theta/p)q}\left(\int_{\Omega_2}\fint_{B(x,2^{-i}) }|\wtil Ef(x)-\wtil Ef(y)|^pd\mu(y)d\mu(x)\right)^{q/p}\\
    &\qquad+\sum_{i\in\N}2^{i(\alpha +\theta/p)q}\left(\int_{\Omega_2}|\wtil Ef(x)|^p\fint_{B(x,2^{-i}) }|\Phi (x)-\Phi (y)|^pd\mu(y)d\mu(x)\right)^{q/p}\\
    &\lesssim \|\wtil Ef\|_{HB^{\alpha +\theta /p}_{p,q}(\Omega,\mu)}^q+\sum_{i\in\N}2^{i(\alpha+ \theta/p)q}\left(\int_{\Omega_2}|\wtil Ef(x)|^p\fint_{B(x,2^{-i}) }d(x,y)^p d\mu(y)d\mu(x)\right)^{q/p}\\
    &\leq \|\wtil Ef\|_{HB^{\alpha +\theta /p}_{p,q}(\Omega,\mu)}^q+\sum_{i\in\N}2^{i(\alpha+ \theta/p)q}\left(\int_{\Omega_2} 2^{-ip}|\wtil Ef(x)|^pd\mu(x)\right)^{q/p}\\
    &\leq     \|\wtil Ef\|_{HB^{\alpha +\theta /p}_{p,q}(\Omega,\mu)}^q+\sum_{i\in\N}2^{i(\alpha+\theta/p-1)q}\left(\int_{\Omega_2}|\wtil Ef(x)|^pd\mu(x)\right)^{q/p}\\
     &\lesssim \|f\|_{HB^\alpha_{p,q}(\partial\Omega,\nu)}^q+\|f\|_{L^p(\partial\Omega,\nu)}^q,
\end{align*}
where the last step follows from the assumption that $\alpha+\theta /p<1$ for the convergence of the series, as well as estimate \eqref{eq:ExtL^pBound} and Theorem \ref{thm:Extension}.

In order to estimate $II$, we note that for $i\in\N$, $x\in\Omega\setminus\Omega_2$, and $y\in B(x,2^{-i})\cap \Omega_2$, we have by the triangle inequality that
\begin{align*}
|Ef(x)-Ef(y)|\le|\Phi(x)||\wtil Ef(x)-\wtil Ef(y)|+|\wtil Ef(y)||\Phi(x)-\Phi(y)|=|\wtil Ef(y)||\Phi(x)-\Phi(y)|,
\end{align*}
since $\Phi (x)=0$. By the doubling property of $\mu$, it also follows that for such $x$ and $y$, 
\begin{equation}\label{eq:doubB(x,2^i)}
    \mu( B(y,2^{-i})  )\leq \mu (B(x,2^{-i+1})  )\leq C_\mu \mu (B(x,2^{-i})  ).
\end{equation}
Hence by \eqref{eq:doubB(x,2^i)}, Tonelli's theorem, and the $1$-Lipschitz condition of $\Phi$, we have
\begin{align*}
II&\lesssim\sum_{i\in\N}2^{i(\alpha+\theta/p)q}\left(\int_{\Omega\backslash\Omega_2}\int_{B(x,2^{-i})\cap \Omega_2}\frac{|\wtil Ef(y)|^p|\Phi(x)-\Phi(y)|^p}{\mu (B(y,2^{-i})  )}d\mu(y)d\mu(x)\right)^{q/p}\\
&= \sum_{i\in\N}2^{i(\alpha+\theta/p)q}\left(\int_{\Omega_2}\int_{B(y,2^{-i})\cap(\Omega\backslash\Omega_2)}\frac{|\wtil Ef(y)|^p|\Phi(x)-\Phi(y)|^p}{\mu (B(y,2^{-i})  )}d\mu(x)d\mu(y)\right)^{q/p}\\
&\le \sum_{i\in\N}2^{i(\alpha+\theta/p)q}\left(\int_{\Omega_2}|\wtil Ef(y)|^p\int_{B(y,2^{-i}) }\frac{d(x,y)^p}{\mu (B(y,2^{-i})  )}d\mu(x)d\mu(y)\right)^{q/p}\\
&\lesssim \sum_{i\in\N}2^{i(\alpha+\theta/p-1)q}\left(\int_{\Omega_2}|\wtil Ef(y)|^pd\mu(y)\right)^{q/p} \lesssim \|f\|_{L^p(\partial\Omega,\nu)}^q.
\end{align*}
It then follows by the estimates of $I$, $II$, and \eqref{eq:InhomExt Energy bound} that
 \[ \|Ef\|_{HB^{\alpha+\theta/p}_{p,q}(\Omega,\mu)}\le C\left(\|f\|_{HB^\alpha_{p,q}(\partial\Omega,\nu)}+\|f\|_{L^p(\partial\Omega,\nu)}\right),
 \]
 where $C$ depends on $\alpha$, $p$, $q$, $\theta$, $C_\mu$ and $C_\theta$.  This completes the proof.

 Now let $q=\infty$. In this case, set
 \[
I:=\sup_{i\in\N}2^{i(\alpha+\theta/p)}\left(\int_{\Omega_2}\fint_{B(x,2^{-i}) }|Ef(x)-Ef(y)|^pd\mu(y)d\mu(x)\right)^{1/p}
\]
and
\[
II:=\sup_{i\in\N}2^{i(\alpha+\theta/p)}\left(\int_{\Omega\setminus\Omega_2}\int_{B(x,2^{-i})\cap \Omega_2}\frac{|Ef(x)-Ef(y)|^p}{\mu(B(x,2^{-i}) )}d\mu(y)d\mu(x)\right)^{1/p}.
\]
 By the same triangle inequality arguments as before,
\begin{align*}
    I&\lesssim 
     \|\wtil Ef\|_{HB^{\alpha +\theta /p}_{p,\infty}(\Omega,\mu)}+\sup_{i\in\N}2^{i(\alpha+ \theta/p)}\left(\int_{\Omega}|\wtil Ef(x)|^p\fint_{B(x,2^{-i}) }d(x,y)^p d\mu(y)d\mu(x)\right)^{1/p}\\
     &\lesssim \|f\|_{HB^\alpha_{p,\infty}(\partial\Omega,\nu)}+\sup_{i\in\N}2^{i(\alpha +\theta/p -1)} \|f\|_{L^p(\partial\Omega,\nu)}=  \|f\|_{HB^\alpha_{p,\infty}(\partial\Omega,\nu)}+2^{\alpha +\theta/p -1} \|f\|_{L^p(\partial\Omega,\nu)}.
\end{align*}
 Similar adaptations from the case $1\le q<\infty$ yield the corresponding estimate for $II$.
\end{proof}

Finally, we show that the trace of this extension operator recovers the original function. 

\begin{lem}\label{lem:TraceOfExtension}
Let $f\in B^\alpha_{p,q}(\partial\Omega,\nu)$ and let $Ef$ be defined as in Theorem~\ref{thm:Extension}.  Then for $\nu$-a.e.\ $z\in\partial\Omega$, we have that 
\[
\lim_{r\to 0^+}\fint_{B(z,r) }|Ef-f(z)|d\mu=0.
\]
\end{lem}
\begin{proof}
As $f\in B^\alpha_{p,q}(\partial\Omega,\nu)\subset L^p(\partial\Omega,\nu)$, $f$ is locally integrable on $\partial\Omega$.  Thus, by the Lebesgue differentiation theorem, $\nu$-a.e.\ $z\in\partial\Omega$ is a Lebesgue point of $f$.  Fix such a point $z\in\partial\Omega$ and $0<r<1$. Notice that for such $r$ we always have $Ef(x) =  \wtil Ef (x)$ for each $x\in B(z,r)$.  By the partition of unity, we then have that 
\begin{align*}
\int_{B(z,r) }&|Ef-f(z)|d\mu=\int_{B(z,r) }\Bigg|\sum_{i,j}\left(\fint_{U_{i,j}}fd\nu\right)\phii_{i,j}(x)-f(z)\Bigg|d\mu(x)\\
	&\le\int_{B(z,r) }\sum_{i,j}\left(\fint_{U_{i,j}}|f-f(z)|d\nu\right)\phii_{i,j}(x)d\mu(x)\\
	&\le\sum_{\substack{l,m\st\\B_{l,m}\cap B(z,r)\ne\varnothing}}\int_{B_{l,m} }\sum_{i,j}\left(\fint_{U_{i,j}}|f-f(z)|d\nu\right)\phii_{i,j}(x)d\mu(x)\\
	&=\sum_{\substack{l,m\st\\B_{l,m}\cap B(z,r)\ne\varnothing}}\int_{B_{l,m} }\sum_{\substack{i,j\st\\2B_{i,j}\cap B_{l,m}\ne\varnothing}}\left(\fint_{U_{i,j}}|f-f(z)|d\nu\right)d\mu(x).
\end{align*}
By bounded overlap, there are at most $C$ indices $(i,j)$ such that $2B_{i,j}\cap B_{l,m}\ne\varnothing$, with $C$ depending only on the doubling constant.  Furthermore, for such $(i,j)$, we have that $U_{i,j}\subset U_{l,m}^*$, with $U_{l,m}^*$ defined as in \eqref{eq:U^*}, and so by doubling, it follows that $\nu(U_{i,j})\simeq\nu(U_{l,m}^*)$.  Hence, we have that 
\begin{align*}
\int_{B(z,r) }|Ef-f(z)|d\mu&\lesssim\sum_{\substack{l,m\st\\B_{l,m}\cap B(z,r)\ne\varnothing}}\int_{B_{l,m} }\fint_{U_{l,m}^*}|f-f(z)|d\nu d\mu(x)\\
	&=\sum_{\substack{l,m\st\\B_{l,m}\cap B(z,r)\ne\varnothing}}\frac{\mu(B_{l,m} )}{\nu(U^*_{l,m})}\int_{U_{l,m}^*}|f-f(z)|d\nu\\
	&\lesssim\sum_{\substack{l,m\st\\B_{l,m}\cap B(z,r)\ne\varnothing}}r_{l,m}^\theta\int_{U_{l,m}^*}|f-f(z)|d\nu\le\sum_{\substack{l,m\st\\B_{l,m}\cap B(z,r)\ne\varnothing}}2^{\theta l}\int_{U_{l,m}^*}|f-f(z)|d\nu.
\end{align*}
If $(l,m)$ is such that $B_{l,m}\cap B(z,r)\ne\varnothing$, then 
\[ 
r\ge d(p_{l,m},z)-r_{l,m}\ge d(p_{l,m},\partial\Omega )-r_{l,m}\geq 7r_{l,m}\ge 7(2^{l-1}),
\] 
and so we have that $2^l<r$.  Let $l_r\in\Z$ be such that $2^{l_r-1}\le r<2^{l_r}$.  Also, if $\zeta\in U^*_{l,m}$ for such $(l,m)$, it follows that 
\begin{align*} 
    d(\zeta,z)&\le d(\zeta,q_{l,m})+d(q_{l,m},p_{l,m})+d(p_{l,m},z)\\
    &<32r_{l,m}+8r_{l,m}+r_{l,m}+r\le 42 r,
\end{align*}
and so $U^*_{l,m}\subset B(z,42r)$.
We then have by Lemma~\ref{lem:bounedoverlapUij} that 
\begin{align*}
\int_{B(z,r) }|Ef-f(z)|d\mu&\lesssim\sum_{l=-\infty}^{l_r-1}2^{\theta l}\sum_{\substack{m\st\\B_{l,m}\cap B(z,r)\ne\varnothing}}\int_{U_{l,m}^*}|f-f(z)|d\nu\\
    &\lesssim\sum_{l=-\infty}^{l_r-1}2^{\theta l}\int_{B(z,42r)\cap\partial\Omega}|f-f(z)|d\nu\\
    &\lesssim r^\theta\int_{B(z,42r)\cap\partial\Omega}|f-f(z)|d\nu.
\end{align*}
Therefore by the $\theta$ codimensional Ahlfors regularity and doubling property of $\nu$, and since $z$ is a Lebesgue point of $f$, we have that 
\begin{align*}
\fint_{B(z,r) }|Ef-f(z)|d\mu&\lesssim\frac{r^\theta}{\mu(B(z,r) )}\int_{B(z,42r)\cap\partial\Omega}|f-f(z)|d\nu\\
    &\lesssim\fint_{B(z,42r)\cap\partial\Omega}|f-f(z)|d\nu\to 0
\end{align*}
as $r\to 0^+$.
\end{proof}

 Notice that the above lemma also holds when  $Ef$ is replaced by $  \wtil Ef$, since for sufficiently small $r$, we have that $Ef(x)=  \wtil Ef(x)$ for all $x\in B(z,r)$. Likewise the following result also holds for $\wtil E$, although we state it for $E$. In both cases, if $  \wtil E$ is considered, then we can consider $f$ to be in $HB_{p,q}^\alpha (\Omega ,\mu )$.

\begin{cor}\label{cor:Identity}
Under the assumptions of Theorem~\ref{thm:Trace}, we have that $T\circ E=\Id$.  That is, for each $f\in B^{\alpha}_{p,q}(\partial\Omega,\nu)$, we have that $TEf(z)=f(z)$ for $\nu$-a.e.\ $z\in\partial\Omega$.
\end{cor}

\begin{proof}
    By \eqref{eq:TraceCondition} and Lemma~\ref{lem:TraceOfExtension}, we have that 
    \[
    \lim_{r\to 0^+}\fint_{B(z,r) }|Ef-TEf(z)|d\mu=0=\lim_{r\to 0^+}\fint_{B(z,r) }|Ef-f(z)|d\mu
    \]
    for $\nu$-a.e.\ $z\in\partial\Omega$.  Therefore, it follows that $TEf(z)=f(z)$ for $\nu$-a.e.\ $z\in\partial\Omega$.
\end{proof}

\section{Hyperbolic fillings and the proof of Theorem~\ref{thm:Intro gen trace}}\label{sec:gen trace thm}

In Section~\ref{section:Trace}, we established trace theorems under the geometric assumption that $\Omega$ is a uniform domain in its completion. In this section, we establish the corresponding result, Theorem~\ref{thm:Intro gen trace}, without this assumption. To do so, we consider uniformized hyperbolic fillings of a metric measure space, as constructed in \cite{BBS}.

\subsection{Uniformized hyperbolic fillings}\label{sec:hyp fill}

Let $(Z,d,\nu)$ be a compact metric measure space, with $\nu$ a doubling measure.  By rescaling the metric if necessary, we may assume that $\diam(Z)<1$.  We will realize $(Z,d,\nu)$, up to a bi-Lipschitz change in the metric, as the boundary of a uniform domain, which is a metric graph obtained through the following construction. 

Fix $z_0\in Z$, and let $A_0:=\{z_0\}$.  For each $n\in\N$, let $A_n\subset Z$ be a maximal $2^{-n}$-separated subset of $Z$, chosen so that $A_n\subset A_{n+1}$.  We obtain the vertex set $V$ of the desired graph by associating to each point in $A_n$ its corresponding level.  That is, 
\[
V:=\bigcup_{n=0}^\infty\bigcup_{z\in A_{n}}(z,n).
\]
The edge relationship $\sim$ between vertices is defined as follows.  Two vertices $(z,n)$ and $(y,m)$ are joined by an edge, i.e.\ $(z,n)\sim (y,m)$, if and only if $n=m$ and $B(z, 2^{-n+1})\cap B(y,2^{-n+1})\ne\varnothing$ or $m=n\pm 1$ and $B(z,2^{-n})\cap B(z,2^{-m})\ne\varnothing$.  We obtain a metric graph $X$ by associating a copy of the unit interval to each edge, and we equip $X$ with the path metric $d_X$.  We denote by $[v,w]$ the interval associated to the edge joining vertices $v$ and $w$.  The metric space $(X,d_X)$ is the \emph{hyperbolic filling} of $(Z,d)$.

Letting $\eps:=\log 2$ and $v_0:=(z_0,0)$, we define the uniformized metric $d_\eps$ on $X$ by
\[
d_\eps(x,y):=\inf_{\gamma}\int_{\gamma}e^{-\eps d_X(\cdot,v_0)}ds
\]
where the infimum is over all paths $\gamma$ in $X$ with endpoints $x$ and $y$.  We denote by $\overline X_\eps$ the completion of $X$ with respect to the metric $d_\eps$, and we denote $\partial_\eps X:=\overline X_\eps\setminus X$.  With this choice of $\eps$, $(\overline X_\eps,d_\eps)$ is bounded. 

For each $\beta>0$, we also lift up the measure $\nu$ to essentially a weighted measure $\mu_\beta$ on $X$ as follows.  We first define the vertex weights 
\[
\hat\mu_\beta((z,n)):=e^{-\beta n}\nu(B(z,2^{-n})),
\]
and then for each $A\subset X$, we define 
\[
\mu_\beta(A):=\sum_{v\in V}\sum_{w\sim v}(\hat\mu_\beta(v)+\hat\mu_\beta(w))\Ha^1([v,w]\cap A).
\]
The following properties of $(\overline X_\eps,d_\eps,\mu_\beta)$ were obtained in \cite[Propositions 4.4, 4.5, 4.6 and Theorem 10.3]{BBS}:

\begin{thm}\label{thm:hyp fill properties}
    Let $(Z,d,\nu)$ be a compact metric measure space, with $\nu$ a doubling measure, and let $\eps:=\log 2$.  Then for each $\beta>0$, the metric measure space $(\overline X_\eps,d_\eps,\mu_\beta)$ satisfies the following:
    \begin{enumerate}
        \item[(i)] $(Z,d)$ is bi-Lipschitz equivalent with $(\partial_\eps X,d_\eps)$,
        \item[(ii)] $(\overline X_\eps,d_\eps)$ is geodesic,
        \item[(iii)] $X$ is a length uniform domain in $\overline X_\eps$, see Remark~\ref{rem:Uniform Domain},
        \item[(iv)] $\mu_\beta$ is doubling,
        \item[(v)] For all $z\in Z$, $0<r<2\diam(Z)$, we have 
        \[
        \nu(B(z,r)\cap Z)\simeq\frac{\mu_\beta(B(z,r))}{r^{\beta/\eps}}.
        \]
    \end{enumerate}
    In each of the above, the comparison constants depend only on $\eps$, $\beta$, and the doubling constant of $\nu$. 
\end{thm}

In a slight abuse of notation, we do not distinguish between balls defined with respect to $d$ and $d_\eps$ due to the bi-Lipschitz equivalence between $(Z,d)$ and $(\partial_\eps X,d_\eps)$. Furthermore, for $A\subset\overline X_\eps$, we interpret $\nu(A)$ to mean $\nu(A\cap Z)$, for ease of notation. 

In \cite{BBS}, it was further shown that $(\overline X_\eps,d_\eps,\mu_\beta)$ supports a $1$-Poincar\'e inequality and that bounded linear trace and extension operators exist between Besov spaces on $Z$ and Newton-Sobolev spaces on $\overline X_\eps$.  As these properties are not relevant to us, we have not included them in the statement of the above theorem.

\subsection{Trace theorem}
We now prove Theorem~\ref{thm:Intro gen trace}, restated below as Theorem~\ref{thm:gen trace}. To do so, we will use the following lemma from \cite{GKS1}:

\begin{lem}[\cite{GKS1}, Lemma 3.10]\label{lem:Hausdorff null} Let $(X,d,\mu)$ be a doubling metric measure space.  Let $f\in L^p(X)$, $0<t<p$, and $M>0$.  Then $\Ha^{-t}(E_M)=0$, where 
\[
E_M=\left\{x\in U:\limsup_{r\to 0^+}r^t\fint_{B(x,r)}|f|^pd\mu>M^p\right\}.
\]
    
\end{lem}

\begin{thm}\label{thm:gen trace}
    Let $(Z,d,\nu)$ be a locally compact, non-complete, bounded metric measure space, with $\nu$ a doubling measure.  Suppose that $\partial Z:=\overline Z\setminus Z$, the boundary of $Z$, is equipped with a Borel measure $\pi$, which is codimension $\theta$ Ahlfors regular with respect to $\nu$ for some $\theta >0$.  Let $1\le p<\infty$, $1\le q\le\infty$, and $\theta/p<\alpha<1$.  Then there exists a bounded, linear trace operator
    \[
    T:B^\alpha_{p,q}(Z,\nu)\to B^{\alpha-\theta/p}_{p,q}(\partial Z,\pi),
    \]
    such that for all $u\in B^\alpha_{p,q}(Z,\nu)$, we have
    \begin{equation}\label{eq:general trace condition}
        \lim_{r\to 0^+}\fint_{B(z,r)}|u-Tu(z)|^pd\nu=0
    \end{equation}
    for $\pi$-a.e.\ $z\in \partial Z$, and
    
    \begin{equation}\label{eq:gen trace energy and Lp}
        \|Tu\|_{HB^{\alpha-\theta/p}_{p,q}(\partial Z,\pi)}\le C\|u\|_{HB^\alpha_{p,q}(Z,\nu)},\qquad\|Tu\|_{L^p(\partial Z,\pi)}\le C\left(\|u\|_{L^p(Z,\nu)}+\|u\|_{HB^\alpha_{p,q}(Z,\nu)}\right),
    \end{equation}
    where $C\ge 1$ depends only on $\alpha$, $p$, $q$, $\theta$, $C_\nu$, $C_\theta$, and $\diam(Z)$.

     Furthermore, the bounded linear extension operator $\wtil E:B^{\alpha-\theta/p}_{p,q}(\partial Z,\pi)\to B^\alpha_{p,q}(Z,\nu)$, given by Corollary~\ref{cor:inhom ext bounded domain}, is a right-inverse of $T$.  That is, for all $f\in B^{\alpha-\theta/p}_{p,q}(\partial Z,\pi)$, we have that 
    \begin{equation}\label{eq:gen TEf=f}
    T\wtil Ef(z)=f(z)
    \end{equation}
    for $\pi$-a.e.\ $z\in\partial Z$.
    
\end{thm}

\begin{proof} 
Consider the uniformized hyperbolic filling $(\overline X_\eps,d_\eps,\mu_\beta)$ of $(\overline Z,d,\nu)$, constructed with parameters $\eps=\log 2$ and $\beta>0$ chosen so that $\sigma:=\beta/\eps=p(1-\alpha)/2$, hence $\sigma<p(1-\alpha)$.  Note that by Theorem~\ref{thm:hyp fill properties}, $\mu:=\mu_\beta$ is doubling, and the measure $\nu$ is codimension $\sigma$ Ahlfors regular with respect to $\mu$, and $(\overline X_\eps,d_\eps)$ is bounded.  By Corollary~\ref{cor:inhom ext bounded domain}, there exists a bounded, linear extension operator
\[
E_Z:B^\alpha_{p,q}(\overline Z,\nu)\to B^{\alpha+\sigma/p}_{p,q}(\overline X_\eps,\mu).
\]
Note that the codimensionality between $\mu$ and $\nu$ implies that $\mu(\overline Z)=0$, and so the above holds since $B^{\alpha+\sigma/p}_{p,q}(\overline X_\eps,\mu)=B^{\alpha+\sigma/p}_{p,q}(X_\eps,\mu)$.

Since $X$ is a length uniform domain in $\overline X_\eps$, the corresponding uniform curves can be extended to the boundary $\overline Z$ by the Arzel\`a-Ascoli theorem, see Remark~\ref{rem:Uniform Domain}.  Thus, it follows that $\overline X_\eps\setminus\partial Z$ is an $A$-uniform domain in $\overline X_\eps$, with $A$ depending only on $\alpha$, $p$, and $C_\nu$.  Furthermore, the measure $\pi$ is codimension $(\sigma+\theta)$ Ahlfors regular with respect to $\mu_\beta$.  Thus, by Theorem~\ref{thm:Intro Inhom Trace}, there exists a bounded, linear trace operator 
\[
T_{\partial Z}:B^{\alpha+\sigma/p}_{p,q}(\overline X_\eps,\mu)\to B^{\alpha-\theta/p}_{p,q}(\partial Z,\pi),
\]
again using the fact that $\mu(\partial Z)=0$.
Define $T:=T_{\partial Z}\circ E_Z$, and let $u\in B^\alpha_{p,q}(Z,\nu)$. The estimates \eqref{eq:gen trace energy and Lp} follow immediately from \eqref{eq:hom trace energy bound}, \eqref{eq:trace inhom Lp}, \eqref{eq:hom extension energy bound}, and \eqref{eq:inhom Lp bounds}.

We now show \eqref{eq:general trace condition}.  To this end, let $0<r<1$.  By the definition of $T$ and Proposition~\ref{prop:tracelimit}, it follows that for $\pi$-a.e.\ $z\in \partial Z$,
\begin{equation}\label{eq:Tu(z)}
Tu(z)=\lim_{\rho\to 0^+}\fint_{B(z,\rho)}E_Zu\,d\mu.
\end{equation}
Fix such a $z\in \partial Z$.  By Theorem~\ref{thm:Trace} and Corollary~\ref{cor:Identity}, there exists a bounded, linear trace operator 
\begin{equation*}\
T_Z:HB^{\alpha+\sigma/p}_{p,q}(\overline X_\eps,\mu)\to HB^\alpha_{p,q}(\overline Z,\nu).
\end{equation*}
such that $T_Z\circ E_Z=\Id$. Again using Proposition~\ref{prop:tracelimit}, it follows that 
\begin{equation}\label{eq:u(w)}
u(w)=\lim_{\rho\to 0^+}\fint_{B(w,\rho)}E_Z u\,d\mu
\end{equation}
for $\nu$-a.e.\ $w\in B(z,r)\cap \overline Z$.  As $\overline X_\eps\setminus \overline Z$ is an $A$-uniform domain, with $A$ depending only on $\alpha$, $p$, and $C_\nu$, we can join $z$ and $w$ by a chain of balls $\{B_k:=B(x_k,r_k)\}_{k\in\Z}$ in  $\overline X_\eps$, as described in Lemma~\ref{lem:uniformcover}.  Similar to the proof of Proposition~\ref{thm:TraceEnergy}, we have by \eqref{eq:Tu(z)}, \eqref{eq:u(w)}, and the properties of this chain of balls that
\begin{align*}
    |u(w)-Tu(z)|^p&\le\left(\sum_{k\ge 0}|(E_Zu)_{B_{k+1}}-(E_Zu)_{B_k}|+\sum_{k<0}|(E_Zu)_{B_{k+1}}-(E_Zu)_{B_k}|\right)^p\\
    &\lesssim\left(\sum_{k\ge 0}\fint_{2B_k}\fint_{2B_k}|E_Zu(x)-E_Zu(y)|d\mu(y)d\mu(x)\right)^p\\
    &\qquad+\left(\sum_{k< 0}\fint_{2B_k}\fint_{2B_k}|E_Zu(x)-E_Zu(y)|d\mu(y)d\mu(x)\right)^p.
\end{align*}

Set $\delta:=\theta/p$.  As in the proof of Proposition~\ref{thm:TraceEnergy} (with $\delta$ in place of $\beta$ there), we then have 
\begin{align*}
    \Bigg(\sum_{k\ge 0}\fint_{2B_k}\fint_{2B_k}&|E_Zu(x)-E_Zu(y)|d\mu(y)d\mu(x)\Bigg)^p\\
    &\lesssim d(z,w)^{\delta p}\int_{C_{z,w}^1}\int_{B(x,d(x,\overline Z))}\frac{|E_Zu(x)-E_Zu(y)|^p}{d(x,\overline Z)^{\delta p}\mu(B(x,d(x,\overline Z)))^2}d\mu(y)d\mu(x)
\end{align*}
and 
\begin{align*}
    \Bigg(\sum_{k< 0}\fint_{2B_k}\fint_{2B_k}&|E_Zu(x)-E_Zu(y)|d\mu(y)d\mu(x)\Bigg)^p\\
    &\lesssim d(z,w)^{\delta p}\int_{C_{z,w}^2}\int_{B(x,d(x,\overline Z))}\frac{|E_Zu(x)-E_Zu(y)|^p}{d(x,\overline Z)^{\delta p}\mu(B(x,d(x,\overline Z)))^2}d\mu(y)d\mu(x),
\end{align*}
 where $C_{z,w}^1:=\bigcup_{k\ge 0}2B_k$ and $C_{z,w}^2:=\bigcup_{k<0}2B_k$.  Therefore, it follows that 
\begin{align}\label{eq:General Trace I and II}
    \fint_{B(z,r)}|u&-Tu(z)|^pd\nu\nonumber\\
    &\lesssim\fint_{B(z,r)}d(z,w)^{\delta p}\int_{C_{z,w}^1}\int_{B(x,d(x,\overline Z))}\frac{|E_Zu(x)-E_Zu(y)|^p}{d(x,\overline Z)^{\delta p}\mu(B(x,d(x,\overline Z)))^2}d\mu(y)d\mu(x)d\nu(w)\nonumber\\
    &+\fint_{B(z,r)}d(z,w)^{\delta p}\int_{C_{z,w}^2}\int_{B(x,d(x,\overline Z))}\frac{|E_Zu(x)-E_Zu(y)|^p}{d(x,\overline Z)^{\delta p}\mu(B(x,d(x,\overline Z)))^2}d\mu(y)d\mu(x)d\nu(w)\nonumber\\
    &=:I+II.
\end{align}
Using the $\sigma$-codimensionality between $\nu$ and $\mu$, as well as Tonelli's theorem, we have
\begin{align*}
    I&\lesssim\frac{r^\sigma}{\mu(B(z,r))}\int_{B(z,r)}\int_{C_{z,w}^1}\int_{B(x,d(x,\overline Z))}\frac{|E_Zu(x)-E_Zu(y)|^pd(z,w)^{\delta p}}{d(x,\overline Z)^{\delta p}\mu(B(x,d(x,\overline Z)))^2}d\mu(y)d\mu(x)d\nu(w)\\
    &\lesssim r^\sigma\fint_{B(z,Cr)}\int_{B(x,d(x,\overline Z))}\frac{|E_Zu(x)-E_Zu(y)|^p}{d(x,\overline Z)^{\delta p}\mu(B(x,d(x,\overline Z)))^2}\int_{B(z,r)}d(z,w)^{\delta p}\chi_{C_{z,w}^1}(x)d\nu(w)d\mu(y)d\mu(x).
\end{align*}
For each $i\in\Z$, let $X_i:=\{x\in\overline X_\eps:2^{i-1}\le d(x,\overline Z)<2^{i}\}$.  Recalling that $0<r<1$, we then have 
\begin{align*}
    I&\lesssim r^{\sigma}\fint_{B(z,Cr)}\sum_{i=-\infty}^{N_C}\chi_{X_i}(x)\fint_{B(x,2^i)}\frac{|E_Zu(x)-E_Zu(y)|^p}{2^{i\delta p}\mu(B(x,2^i))}\int_{B(z,r)}d(z,w)^{\delta p}\chi_{C_{z,w}^1}(x)d\nu(w)d\mu(y)d\mu(x)\\
    &\lesssim r^{\sigma+\delta p}\fint_{B(z,Cr)}\sum_{i=-\infty}^{N_C}2^{-i\delta p}\fint_{B(x,2^i)}|E_Zu(x)-E_Zu(y)|^p\frac{\nu(B(x,2^i)\cap \overline Z)}{\mu(B(x,2^i))}d\mu(y)d\mu(x)\\
    &\lesssim r^{\sigma+\delta p}\fint_{B(z,Cr)}\sum_{i=-\infty}^{N_C}2^{-i(\delta p+\sigma)}\fint_{B(x,2^i)}|E_Zu(x)-E_Zu(y)|^pd\mu(y)d\mu(x).
\end{align*}
Similarly, we obtain
\begin{align*}
    II\lesssim r^{\sigma+\delta p}\fint_{B(z,Cr)}\sum_{i=-\infty}^{N_C}2^{-i(\delta p+\sigma)}\fint_{B(x,2^i)}|E_Zu(x)-E_Zu(y)|^pd\mu(y)d\mu(x),
\end{align*}
By our choice of $\delta$ and \eqref{eq:General Trace I and II}, we then have 
\begin{align}\label{eq:Haus null ineq}
    \fint_{B(z,r)}|u-Tu(z)|^pd\nu&\lesssim r^{\sigma+\theta}\fint_{B(z,Cr)}
    \sum_{i=-\infty}^{N_C}2^{-i(\sigma+\theta)}\fint_{B(x,2^i)}|E_Zu(x)-E_Zu(y)|^pd\mu(y)d\mu(x)\nonumber\\
    &=:r^{\sigma+\theta}\fint_{B(z,Cr)}g(x)^pd\mu(x),
\end{align}
where we have defined the function 
\[
g(x):=\left(\sum_{i=-\infty}^{N_C}2^{-i(\sigma+\theta)}\fint_{B(x,2^i)}|E_Zu(x)-E_Zu(y)|^pd\mu(y)\right)^{1/p}.
\]

We claim that $g\in L^p(\overline X_\eps,\mu)$.  To show this, let $\tau:=\alpha-\theta/p$ and note that $\tau>0$ by our assumptions.  We then have by H\"older's inequality and the boundedness of the operator $E_Z$ that 
\begin{align*}
    \|g\|_{L^p(\overline X_\eps,\mu)}&=\left(\int_{\overline X_\eps}\sum_{i=-\infty}^{N_C}2^{-i(\sigma+\theta)}\fint_{B(x,2^i)}|E_Zu(x)-E_Zu(y)|^pd\mu(y)d\mu(x)\right)^{1/p}\\
    &\le\sum_{i=-\infty}^{N_C}2^{-i(\sigma+\theta)/p}\left(\int_{\overline X_\eps}\fint_{B(x,2^i)}|E_Zu(x)-E_Zu(y)|^pd\mu(y)d\mu(x)\right)^{1/p}\\
    &\le\sum_{i=-\infty}^{N_C}2^{-i((\sigma+\theta)/p+\tau)}2^{i\tau}\left(\int_{\overline X_\eps}\fint_{B(x,2^i)}|E_Zu(x)-E_Zu(y)|^pd\mu(y)d\mu(x)\right)^{1/p}\\
    &\lesssim\left(\sum_{i=-\infty}^{N_C}2^{-i((\sigma+\theta)/p+\tau)q}\left(\int_{\overline X_\eps}\fint_{B(x,2^i)}|E_Zu(x)-E_Zu(y)|^pd\mu(y)d\mu(x)\right)^{q/p}\right)^{1/q}\\
    &=\left(\sum_{i=-\infty}^{N_C}2^{-i(\alpha+\sigma/p)q}\left(\int_{\overline X_\eps}\fint_{B(x,2^i)}|E_Zu(x)-E_Zu(y)|^pd\mu(y)d\mu(x)\right)^{q/p}\right)^{1/q}\\
    &\le\|E_Zu\|_{HB^{\alpha+\sigma/p}_{p,q}(\overline X_\eps,\mu)}<\infty.
\end{align*}
Though the above estimate is written for $1\le q<\infty$, it holds with the standard modifications for $q=\infty$.

Since $\theta<\alpha p$ and $\sigma<p(1-\alpha)$, we have that $\sigma+\theta<p$. Hence, from Lemma~\ref{lem:Hausdorff null} and \eqref{eq:Haus null ineq}, it follows that 
\[
\lim_{r\to 0^+}\fint_{B(z,r)}|u-Tu(z)|^pd\nu=0
\]
for $\Ha^{-(\sigma+\theta)}_\mu$-a.e.\ $z\in\partial Z$.  As $\pi$ is codimension $(\sigma+\theta)$ Ahlfors regular with respect to $\mu$, we have by Lemma \ref{lem:codimhausdorff} that $\pi\simeq\Ha^{-(\sigma+\theta)}_\mu|_{\partial Z}$, which completes the proof of \eqref{eq:general trace condition}.

It remains to show \eqref{eq:gen TEf=f}.  To this end, let $f\in B^{\alpha-\theta/p}_{p,q}(\partial Z,\pi)$.  Then from Corollary~\ref{cor:inhom ext bounded domain}, it follows that $\wtil Ef\in B^\alpha_{p,q}(\overline Z,\nu)$, and so from \eqref{eq:general trace condition}, we have that for $\pi$-a.e.\ $z\in\partial Z$, 
\begin{equation*}
    \lim_{r\to 0^+}\fint_{B(z,r)}|\wtil Ef-T\wtil Ef(z)|^pd\nu=0.
\end{equation*}
Moreover, from Lemma~\ref{lem:TraceOfExtension}, with $\nu$ and $\pi$ playing the role of $\mu$ and $\nu$ there, respectively, it follows that for $\pi$-a.e.\ $z\in\partial Z$, 
\[
\lim_{r\to 0^+}\fint_{B(z,r)}|\wtil Ef-f(z)|d\nu=0.
\]
By these limits and the triangle inequality, we have that for $\pi$-a.e.\ $z\in\partial Z$, 
\begin{align*}
    |f(z)-T\wtil Ef(z)|\le\lim_{r\to 0^+}\fint_{B(z,r)}|\wtil Ef-T\wtil Ef(z)|^pd\nu+\lim_{r\to 0^+}\fint_{B(z,r)}|\wtil Ef-f(z)|d\nu=0,
\end{align*}
which gives us \eqref{eq:gen TEf=f}.
\end{proof}

\begin{remark}
   As mentioned in the introduction, it may be possible to remove the boundedness assumption in Theorem~\ref{thm:gen trace} by using the hyperbolic filling construction given in \cite{Butler1,Butler2} instead of the construction from \cite{BBS}.  However, as these papers are as-of-yet unpublished, we leave this investigation for future work. 
\end{remark}

\vspace{10pt}
{\fontsize{9pt}{9pt}\selectfont

\noindent I.C.: Institute of Mathematics of the Polish Academy of Sciences, Jędrzeja Śniadeckich 8, 00-656 Warsaw, Poland.

\noindent E-mail: {\tt icaamanoaldemunde@impan.pl}
\vskip.2cm
\noindent J.K.: Dept. of Mathematical Sciences, P.O.~Box 210025, University of Cincinnati, Cincinnati, OH~45221-0025, U.S.A.

\noindent E-mail: {\tt klinejp@ucmail.uc.edu}

}


\begin{thebibliography}{A}
\frenchspacing
%

\bibitem{BBCK}M. Barlow, R. Bass, Z.-Q. Chen, M. Kassmann \emph{Non-local Dirichlet forms and symmetric jump processes.} Trans. Amer. Math. Soc. 361, No. 4, (2009), 1963--1999.

\bibitem{BB}A. Bj\"orn and J. Bj\"orn,
\emph{Nonlinear potential theory on metric spaces}.
EMS Tracts in Mathematics, 17. European Mathematical Society (EMS), Z\"urich, (2011), xii+403 pp.

\bibitem{BBS} A. Björn, J. Björn, N. Shanmugalingam. \textit{Extension and trace results for doubling metric measure spaces and their hyperbolic fillings.} J. Math. Pures Appl.\textbf{159}, (2022), 196--249.

\bibitem{BS} J. Björn, N. Shanmugalingam. \emph{Poincaré inequalities, uniform domains and extension properties for Newton-Sobolev functions in metric spaces.} J. Math. Anal. Appl. \textbf{332}, (2007), no. 1, 190--208.

\bibitem{Butler2} C. Butler. \emph{Extension and trace theorems for noncompact doubling spaces.} Preprint. (2021), https://arxiv.org/pdf/2009.10168

\bibitem{Butler1} C. Butler. \emph{Uniformizing Gromov hyperbolic spaces with Busemann functions.} Preprint. (2021), https://arxiv.org/pdf/2007.11143

\bibitem{CRS} L. Caffarelli, J.M. Roquejoffre, and O. Savin, \emph{Non-local minimal surfaces}. Comm. Pure Appl. Math. \textbf{63}, (2010), 1111--1144.

\bibitem{CS} L. Caffarelli, L. Silvestre. \emph{An extension problem related to the fractional Laplacian.} Comm. Partial Differential
Equations \textbf{32}, (2007), no. 7-9, 1245--1260.

\bibitem{CSt} L. Caffarelli, P. R. Stinga. \emph{Fractional elliptic equations, Caccioppoli estimates and regularity.} Ann. Inst. H.
Poincaré C Anal. Non Linéaire 33 (2016), no. 3, 767–807. 

\bibitem{CGKS} L. Capogna, R. Gibara, R. Korte, N. Shanmugalingam. \emph{Fractional $p$-Laplacians via Neumann problems in unbounded metric measure spaces.} Preprint. (2024), 
https://doi.org/10.48550/arXiv.2410.18883

\bibitem{CGKS2} L. Capogna, R. Gibara, R. Korte, N. Shanmugalingam. \emph{Sharp Hölder regularity of weak solutions of the Neumann problem and applications to nonlocal PDE in metric measure spaces.} Preprint. (2025), 
https://doi.org/10.48550/arXiv.2505.14950
 


\bibitem{CKKSS} L. Capogna, J. Kline, R. Korte, N. Shanmugalingam, M. Snipes, \emph{Neumann problems for $p$-harmonic functions, and induced nonlocal operators in metric measure spaces.} (To appear) Amer. J. Math, (2025).

\bibitem{AChen} A. Chen. \emph{Boundary Harnack Principle on Uniform Domains.} Potential Anal \textbf{62}, (2025), 739--759. https://doi.org/10.1007/s11118-024-10154-4

 \bibitem{C} Z.-Q. Chen. \emph{Multidimensional symmetric stable processes.} Korean J. Comput. Appl. Math. \textbf{6}, (1999), no. 2,
227--266.

\bibitem{D}B. S. Dyda. \emph{A fractional order Hardy inequality.} Illinois J. Math. \textbf{48}, (2004), no. 2, 575--588.

\bibitem{EGKSS} S. Eriksson-Bique, G. Giovannardi, R. Korte, N. Shanmugalingam, G. Speight. \emph{Regularity of solutions to the fractional Cheeger-Laplacian on domains in metric spaces of bounded geometry,} J. Diff. Eq., \textbf{306}, (2022), 590--632.


\bibitem{G} E. Gagliardo. \emph{Caratterizzazioni delle tracce sulla frontiera relative ad alcune classi di funzioni in n variabili.}
Rend. Semin. Mat. Univ. Padova \textbf{27}, (1957), 284--305.

\bibitem{GKS1} R. Gibara, R. Korte, N. Shanmugalingam. \textit{Solving a Dirichlet problem on unbounded domains via a conformal transformation.} Mathematische Annalen,  \textbf{389(3)}, (2024), 2857--2901.

\bibitem{Nages-Ryan}
R. Gibara, N. Shanmugalingam. \emph{Trace and extension theorems for homogeneous Sobolev and Besov spaces for unbounded uniform domains in metric measure spaces}.  Proc. Steklov Inst. Math. \textbf{323}, (2023), 101--119.

\bibitem{GKS} A. Gogatishvili, P. Koskela, N. Shanmugalingam. \textit{Interpolation properties of Besov spaces defined on metric spaces}. Math. Nachr., \textbf{283}, (2010), 215--231.

\bibitem{Gromov} M. Gromov. \textit{Metric Structures for Riemannian and Non-Riemannian Spaces.} Modern Birkhäuser Classics, Birkhäuser Boston, MA, 2007.

\bibitem{HIT} T. Heikkinen, L. Ihnatsyeva, H. Tuominen. \emph{Measure density and extension of Besov and Triebel-Lizorkin functions.} J. Fourier Anal. Appl.  \textbf{22}, (2016), 334--382.

\bibitem{H} J. Heinonen. \emph{Lectures on analysis on metric spaces.} Universitext, Springer, New York, 2001.

\bibitem{HKST}
J. Heinonen, P. Koskela, N. Shanmugalingam, J. T. Tyson,
\emph{Sobolev Spaces on Metric Measure Spaces. An Approach Based on Upper Gradients.}
New Math. Monogr. 27. Cambridge University Press,  Cambridge, 2005.

\bibitem{HiKu} M. Hino, T. Kumagai. \emph{A trace theorem for Dirichlet forms on fractals.} J. Funct. Anal. \textbf{238}, (2006), no. 2, 578--611.

\bibitem{IMV} L. Ihnatsyeva, K. Mohanta, A.V. Vähäkangas. \emph{ Fractional Hardy inequalities and capacity density.} Calc. Var. Partial Differential Equations, Volume 64, \textbf{136}, (2025), doi: 10.1007/s00526-025-02999-3.

\bibitem{JW}
A. Jonsson, H. Wallin, \emph{Function spaces on subsets of $\R^n$} Math. Rep. 2(1), Harwood Academic Publishers, London, 1984. 

\bibitem{KajMur} N. Kajino, M. Murugan. \emph{Heat kernel estimates for boundary traces of reflected diffusions on uniform domains.} Preprint. (2023), https://arxiv.org/abs/2312.08546

\bibitem{KLS} J. Kline, F. Li, N. Shanmugalingam. \textit{Well-posedness of Dirichlet boundary value problems for reflected fractional $p$-Laplace-type inhomogeneous equations in compact doubling metric measure spaces}. Preprint. (2025), https://arxiv.org/abs/2408.02624

\bibitem{LS} P. Lahti, N. Shanmugalingam. \emph{Trace theorems for functions of bounded variation in metric spaces.} J. Funct. Anal. {\bf 274}, (2018), 2431--2460.

\bibitem{M} L. Mal\'y. \emph{Trace and extension theorems for Sobolev-type functions in metric spaces.} Preprint. (2017), https://arxiv.org/abs/1704.06344

\bibitem{MSS} L. Mal\'y, M. Snipes, N. Shanmugalingam.
\emph{Trace and extension theorems for functions of bounded variation.} 
Annali SNS {\bf 18}, (2018), 313--341.

\bibitem{Marcos}
M. A. Marcos. \emph{A trace theorem for Besov functions on spaces of homogeneous type.} Publicacions Matemàtiques, \textbf{62}, (2023), no. 1, 185--211.

\bibitem{Mur} M. Murugan. \emph{Heat kernel for reflected diffusion and extension property on uniform domains} Preprint. (2024), https://arxiv.org/abs/2304.03908

\bibitem{SS} E. Saksman, T. Soto. \emph{Traces of Besov, Triebel-Lizorkin and Sobolev Spaces on Metric Measure Spaces.} Anal. Geom. Metr. Spaces \textbf{5}, (2017), 98--115.

\bibitem{S} N. Shanmugalingam. \emph{Newtonian spaces: an extension of Sobolev spaces to metric measure spaces.} Rev. Mat.
Iberoamericana \textbf{16}, (2000), no. 2, 243--279.



\end{thebibliography}
\end{document}